# ON THE FUNDAMENTAL PROPERTIES OF LINEAR PARAMETER-VARYING DYNAMIC SYSTEMS UNDER PARAMETRICAL MULTI-PERTURBATIONS. APPLICATIONS TO TIME-DELAY SYSTEMS. PRELIMINARY RESULTS

M. De la Sen, Institute for Research and Development of Processes

Department of Electricity and Electronics. **Faculty of Science and Technology**

Campus of Leioa (Bizkaia). Aptdo. 644 de Bilbao

48080-Bilbao. SPAIN

**Abstract**. This paper deals with a unifying approach to the problems of computing the admissible sets of parametrical multi perturbations in appropriate bounded sets such that some fundamental properties of parameter-varying linear dynamic systems are maintained provided that the so-called (i.e. perturbation-free) nominal system possesses such properties. The sets of parametrical multi perturbations include any combinations of parametrical multi perturbations in the matrix of dynamics as well as in the control, output and input-output interconnection matrices which belong to some prescribed bounded domain in the complex space. The various properties which are investigated are controllability, observability, output controllability and existence of minimal state-space realizations together with the associate existence or not of associate decoupling, transmission and invariant zeros. All the matrices of parameters including the nominal and the disturbed ones which parameterize the dynamic system may be real or complex. The radii of the multi-parametrical perturbations are calculated in a simple way. The obtained results are then applied to systems subject to a finite number of point internal delays and parametrical multi perturbations by comparing the state- space descriptions of such systems with the general descriptions previously investigated. In particular, the contributions of the delays to the spectral descriptions are assimilated to the contributions of a set of varying parameters in a domain for the general description.

## 1.     Introduction

The problem of robust stability of dynamic systems has received important attention in the last two decades, [1-3]. The related investigations require in general ad-hoc mathematical tools from Mathematical and Functional Analysis, [1-4]. Recently, the notion of stability radius has been used for related investigations, [5-7]. The stability radius of a linear dynamic system is a positive real number which defines the minimum size, in terms of norm, of a parametrical perturbation, belonging to an admissible class,  such that the resulting system becomes unstable or critically stable provided that the nominal (i.e. perturbation-free) system is stable . Such a characterization has been used successfully in [5-7] to investigate the maximum size of both structured and unstructured multi parametrical (in general, complex) perturbations so that positive  systems are maintained stable provided that its nominal part is stable. Further advantages of focusing the robust stability problem in that way are that the robust stability of wide classes of parameter-varying dynamic systems, including some described by functional equations, may be studied in a unified way, [5-9]. The purpose of this paper is to study the fundamental properties of controllability, observability, stabilizability and detectability of parameter-varying linear systems, in a way inspired in the ideas developed in [5-7]. However, some variations are that the parameter-varying systems under ( in general complex) parametrical multi perturbations are not necessarily positive and that spectral stability radii are not involved since the problem at hand is not that of robust stabilization. The main idea is to maintain the Popov-Belevitch-Hautus matrix functions [1] for the investigated properties being full rank for all multi perturbations provided that the matrix of the nominal system is also full rank. The multi perturbations are of a given structured class on a certain domain where the varying parameters belong to. The worst case of the admissible perturbations establishes the robustness degree of the property. The extension to unstructured perturbations is not discussed since it is direct, and even more simple, than that for structured ones. The study is addressed in a unified way for all the properties. In particular, it has direct interest in realization theory since the size of the disturbances which maintain a minimal state-space realization of the system may be characterized, provided that the nominal realization is minimal. To this end, the best of the two worst cases of losing either controllability or observability by the perturbed system of a nominally controllable and observable system ensures that the state-space realization is still minimal. A direct extension is that if the nominal system is stabilizable and detectable, the best of the two worst cases of losing some of both properties for some multi perturbation in the given class still guarantees that any eventual zero-pole cancellation in the transfer matrix is stable. The technical mechanism employed to investigate the various properties is the construction of square auxiliary matrix functions which are symmetrical positive definite (or Hermitian in the complex case). If their associate minimum singular value becomes zero, or equivalently, if their determinants become zero for some perturbation while their counterparts of the nominal system are positive then the investigated property is lost. One takes advantage that the functions characterizing the singular values and the determinant of a complex continuous matrix function are continuous functions on the definition domain of such a matrix function. The results are easily extendable to parameter-varying dynamic linear internal and external delays under multi – perturbations. There is an important background on time-delay systems, [10-34] including models of neural networks including delays, [25-27]. In particular, the study of stability of time-delay systems has received attention in [10-11], [15], [19-24], [30-34], the positivity and periodicity of the trajectory-solutions have being investigated in [13-14], [27-29], [32], [35-36] and the state-trajectory solutions under impulsive controls in [12-13], [15-17] including the case of singular systems, [13]. Sufficiency–type conditions for the robust characterization of those properties follows directly from the general study based on the fact that time-delay systems might be characterized as nD- systems [8-9], [21]. Refinements of the conditions either allowing to derive results dependent on the delay sizes or guaranteeing that the properties hold for some cases not included in the nD- system characterizations are also investigated. This is addressed by further considering the spectral characterization of linear delay systems under quasi-polynomials.

## 2. Notation

Subsets of the complex and real fields $R$ and $C$ are:

$$R_{+0} := \{r \in R : Re(r) \geq 0\}, \quad R_+ := \{r \in R : Re(r) > 0\}, \quad R_- := \{r \in R : Re(r) < 0\}$$

$$C_{+0} := \{r \in C : Re(r) \geq 0\}, \quad C_+ := \{r \in C : Re(r) > 0\}, \quad C_- := \{r \in C : Re(r) < 0\}$$

which stand for nonnegative, positive and negative real numbers or complex numbers of nonnegative, positive and negative real abscissas, respectively. $C_1 := \{z \in C : |z| = 1\}$ is the unit circumference in the complex plane.

$C^{(q)}(R_{+0}, K^p)$ and $PC^{(q)}(R_{+0}, K^p)$ with $K = C$ or $K = R$ denote, respectively, real or complex vector functions of class $q \geq 0$ and those of class $q-1$ (if $q \geq 1$) whose q-th derivative is continuously differentiable with piecewise continuous $q$-th derivative on the definition domain $R_{+0}$. $PC^{(0)}(R_{+0}, K^p)$ denotes the class of piecewise continuous p-vector functions with domain $R_{+0}$ and range $K^p$.

$I_n$ denotes the n-th identity matrix. The superscript "*" denotes the transpose conjugate of a complex matrix resulting in the usual transpose for real matrices (denoted with the superscript "T"), $det\, M$ and $\sigma(M) := \{\lambda \in C : det(\lambda I_n - M) = 0\}$ denote the determinant of the matrix M and its spectrum, respectively. Subsets of $\sigma(M)$ are $\sigma_+(M) := \{\lambda \in \sigma(M) : Re\,\lambda > 0\}$, $\sigma_{+0}(M) := \{\lambda \in \sigma(M) : Re\,\lambda \geq 0\}$, $\sigma_-(M) := \{\lambda \in \sigma(M) : Re\,\lambda < 0\}$ and $\sigma_{-0}(M) := \{\lambda \in \sigma(M) : Re\,\lambda \leq 0\}$.

The spectral radius and spectral abscissa of M are denoted, respectively, by $\rho(M) := \{max\,|\lambda| : \lambda \in \sigma(M)\}$ and $\mu(M) := \{max\,Re\,\lambda : \lambda \in \sigma(M)\}$. The singular values of M are the positive squares of the eigenvalues of any of the matrix products $MM^*$ and $M^*M$ provided that they exist which are real and satisfy $\bar{\sigma}(M) \geq \underline{\sigma}(M) \geq 0$ where $\bar{\sigma}(M)$ and $\underline{\sigma}(M)$ are the maximum and minimum singular value of M, respectively. Note that at least one of the matrix products $MM^*$ and $M^*M$ always exist. The spectral (or $\ell_2$) vector norm and associated induced matrix norm are denoted by $\|.\|_2$.

The interior, boundary (frontier) and closure of a set $K$ are denoted as $K^0$, $K^{Fr}$, and $cl\,K$, respectively.

$\hat{f}(s) := Lap(f(t))$, $f(t) := Lap^{-1}(\hat{f}(s))$ is a pair of Laplace transform and Laplace anti-transform provided that such an anti-transform exists.

The Kronecker (or direct) product of the matrices $A = (a_{ij})$ and B is denoted by $A \otimes B = (a_{ij} B)$ and for such a matrix A, $vec(A) = (a_1^*, a_2^*, \cdots, a_n^*)^*$ with $a_i^*$ being the i-th row of A, i.e. $A^* = (a_1, a_2, \ldots, a_n)$. The boundary (or frontier) of a set $Q$ is denoted as $Q^{Fr}$.

The complex unit is $i = \sqrt{-1}$.

### 3. The parameter-varying system and associate fundamental properties: controllability, observability, minimal realizability, stabilizability and detectability

The parameter-varying linear time-invariant dynamic system to be considered is

$$\dot{x}(z,t) = \left(A(z^{(A)}) + \tilde{A}(z^{(A)})\right) x(z,t) + \left(B(z^{(B)}) + \tilde{B}(z^{(B)})\right) u(t) \tag{1}$$

$$y(z,t) = \left(C(z^{(C)}) + \tilde{C}(z^{(C)})\right) x(z,t) + \left(D(z^{(D)}) + \tilde{D}(z^{(C)})\right) u(t) \tag{2}$$

, $\forall t \in R_{+0} := R_+ \cup \{0\}$, with $R_+ := \{s \in R : r > 0\}$, subject to initial conditions $x(z,0) = x_0(z) \in C^n$, where $X \subset C^n$, $Y \subset C^p$ and $Y \subset C^m$ are, respectively, the state, output and input linear spaces, and $x \in C^{(1)}(C^q \times R_{+0}, X)$, $y \in C^{(0)}(C^q \times R_{+0}, Y)$ and $u \in PC^{(0)}(C^q \times R_{+0}, U)$ are, the everywhere continuously time-differentiable n-vector state trajectory solution, the piecewise continuous p-vector output-trajectory solution, and the piecewise continuous m-vector control input, respectively, with $p \leq m \leq n$ and $q := q_A + q_B + q_C + q_D + 4$, and

$A: C^{q_A+1} \to C^{n \times n}$, $B: C^{q_B+1} \to C^{n \times m}$, $C: C^{q_C+1} \to C^{p \times n}$ and $D: C^{q_D+1} \to C^{p \times m}$; and

$\tilde{A}: C^{q_A+1} \to C^{n \times n}$, $\tilde{B}: C^{q_B+1} \to C^{n \times m}$, $\tilde{C}: C^{q_C+1} \to C^{p \times n}$ and $\tilde{D}: C^{q_D+1} \to C^{p \times m}$ are,

respectively, the (so-called) nominal and perturbation (complex-valued) matrices of dynamics, control, output and input-output interconnections whose parameter-varying arguments are defined by the respective complex-valued $(q_A+1)$, $(q_B+1)$, $(q_C+1)$ and $(q_D+1)$- tuples:

$$z^{(A)} := \left(1, z_1^{(A)}, z_2^{(A)}, \ldots, z_{q_A}^{(A)}\right) \in \{1\} \times C^{q_A} \subset C^{q_A+1}, \quad z^{(B)} := \left(1, z_1^{(B)}, z_2^{(B)}, \ldots, z_{q_B}^{(B)}\right) \in \{1\} \times C^{q_B} \subset C^{q_B+1}$$

$$\tag{3}$$

$$z^{(C)} := \left(1, z_1^{(C)}, z_2^{(C)}, \ldots, z_{q_C}^{(C)}\right) \in \{1\} \times C^{q_C} \subset C^{q_C+1}, \quad z^{(D)} := \left(1, z_1^{(D)}, z_2^{(D)}, \ldots, z_{q_D}^{(D)}\right) \in \{1\} \times C^{q_D} \subset C^{q_D+1}$$

$$\tag{4}$$

The (so-called) nominal system is affine parameter-varying defined from (1)-(2) with

$$A\left(z^{(A)}\right):=\sum_{i=0}^{q_A} z_i^{(A)} A_i,\ B\left(z^{(B)}\right):=\sum_{i=0}^{q_B} z_i^{(B)} B_i,\ C\left(z^{(C)}\right):=\sum_{i=0}^{q_C} z_i^{(C)} C_i,\ D\left(z^{(D)}\right):=\sum_{i=0}^{q_D} z_i^{(D)} D_i$$

(5)

where $z_0^{(A)}=z_0^{(B)}=z_0^{(C)}=z_0^{(D)}=1$, $\tilde{A}\left(z^{(A)}\right)=0$, $\tilde{B}\left(z^{(B)}\right)=0$, $\tilde{C}\left(z^{(C)}\right)=0$ and $\tilde{D}\left(z^{(D)}\right)=0$, and

$A_i: C \to C^{n\times n}\ (i\in \overline{q}_A^{\,0})$, $B_i: C \to C^{n\times m}\ (i\in \overline{q}_B^{\,0})$, $C_i: C \to C^{p\times n}\ (i\in \overline{q}_C^{\,0})$, $D_i: C \to C^{p\times m}\ (i\in \overline{q}_D^{\,0})$

(6)

with $\overline{n}^{\,0}:=\overline{n}\cup\{0\}$, $\overline{n}:=\{1,2,...,n\}$. The nominal system is the (unperturbed) reference one to then establish the nominal bounded domain (i.e. a connected open set) $C_{\alpha 0}\subset C^q$ where the nominal system fulfills the various investigated properties for all $z:=\left(z^{(A)},z^{(B)},z^{(C)},z^{(D)}\right)\in C_{\alpha 0}$ and also the bounded domain $C_\alpha \subset C_{\alpha 0}$ where a class of systems (1)-(2), eventually submitted to perturbations, still maintain the particular property under investigation kept by the nominal one on $C_{\alpha 0}$. The class of systems (1)-(2), which include the nominal system as particular case, are defined via (3)-(6) for parametrical multi perturbations in a set $P:=P_A\times P_B\times P_C\times P_D\subset C^{q_A+1}\times C^{q_B+1}\times C^{q_C+1}\times C^{q_D+1}\equiv C^q$ of the form:

$$\tilde{A}\left(z^{(A)}\right):=\sum_{i=0}^{q_A} z_i^{(A)}\tilde{A}_i=\sum_{i=0}^{q_A}\sum_{j=1}^{n_A} z_i^{(A)} D_{ij}^{(A)}\Delta_{ij}^{(A)} E^{(A)},\ \tilde{B}\left(z^{(B)}\right):=\sum_{i=0}^{q_B} z_i^{(B)}\tilde{B}_i=\sum_{i=0}^{q_B}\sum_{j=1}^{n_B} z_i^{(B)} D_{ij}^{(B)}\Delta_{ij}^{(B)} E^{(B)}$$

$$\tilde{C}\left(z^{(C)}\right):=\sum_{i=0}^{q_C} z_i^{(C)}\tilde{C}_i=\sum_{i=0}^{q_C}\sum_{j=1}^{n_C} z_i^{(C)} D_{ij}^{(C)}\Delta_{ij}^{(C)} E^{(C)},\ \tilde{D}\left(z^{(D)}\right):=\sum_{i=0}^{q_D} z_i^{(D)}\tilde{D}_i=\sum_{i=0}^{q_D}\sum_{j=1}^{n_D} z_i^{(D)} D_{ij}^{(D)}\Delta_{ij}^{(D)} E^{(D)}$$

(7)

where $\tilde{A}: C^{q_A+1}\to P_A\subset C^{n\times n}$, $\tilde{B}: C^{q_B+1}\to P_B\subset C^{n\times m}$, $\tilde{C}: C^{q_C+1}\to P_C\subset C^{p\times n}$ and $\tilde{D}: C^{q_D+1}\to P_D\subset C^{p\times m}$ with

$\tilde{A}_i: C \to C^{n\times n}\ (i\in \overline{q}_A^{\,0})$, $\tilde{B}_i: C \to C^{n\times m}\ (i\in \overline{q}_B^{\,0})$, $\tilde{C}_i: C \to C^{p\times n}\ (i\in \overline{q}_C^{\,0})$, $\tilde{D}_i: C \to C^{p\times m}\ (i\in \overline{q}_D^{\,0})$

(8)

are defined by $\tilde{A}_i:=\sum_{j=1}^{n_A} D_{ij}^{(A)}\Delta_{ij}^{(A)} E^{(A)}\ (i\in \overline{q}_A^{\,0})$, $\tilde{B}_i:=\sum_{j=1}^{n_B} D_{ij}^{(B)}\Delta_{ij}^{(B)} E^{(B)}\ (i\in \overline{q}_B^{\,0})$, $\tilde{C}_i:=\sum_{j=1}^{n_C} D_{ij}^{(C)}\Delta_{ij}^{(C)} E^{(C)}\ (i\in \overline{q}_C^{\,0})$ and $\tilde{D}_i:=\sum_{j=1}^{n_D} D_{ij}^{(D)}\Delta_{ij}^{(D)} E^{(D)}\ (i\in \overline{q}_A^{\,0})$, where

$D_{ij}^{(A)}: C\to C^{n\times \ell_{ij}^{(A)}},(i,j)\in \overline{q}_A^{\,0}\times \overline{n}_A;\ D_{ij}^{(B)}: C\to C^{n\times \ell_{ij}^{(B)}},(i,j)\in \overline{q}_B^{\,0}\times \overline{n}_B$

$D_{ij}^{(C)}: C\to C^{p\times \ell_{ij}^{(C)}},\forall (i,j)\in \overline{q}_C^{\,0}\times \overline{n}_C;\ D_{ij}^{(D)}: C\to C^{p\times \ell_{ij}^{(D)}},\forall (i,j)\in \overline{q}_D^{\,0}\times \overline{n}_D$

(9)

$$\Delta_{ij}^{(A)}: \boldsymbol{C} \to \boldsymbol{C}^{\ell_{ij}^{(A)} \times \ell^{(A)}}, \forall (i,j) \in \overline{q}_A^0 \times \overline{n}_A; \quad \Delta_{ij}^{(B)}: \boldsymbol{C} \to \boldsymbol{C}^{\ell_{ij}^{(B)} \times \ell^{(B)}}, \forall (i,j) \in \overline{q}_B^0 \times \overline{n}_B$$

$$\Delta_{ij}^{(C)}: \boldsymbol{C} \to \boldsymbol{C}^{\ell_{ij}^{(C)} \times \ell^{(C)}}, \forall (i,j) \in \overline{q}_C^0 \times \overline{n}_C; \quad \Delta_{ij}^{(D)}: \boldsymbol{C} \to \boldsymbol{C}^{\ell_{ij}^{(D)} \times \ell^{(D)}}, \forall (i,j) \in \overline{q}_D^0 \times \overline{n}_D$$

(10)

$$E^{(A)}: \boldsymbol{C} \to \boldsymbol{C}^{\ell_f^{(A)} \times n}; \quad E^{(B)}: \boldsymbol{C} \to \boldsymbol{C}^{\ell_f^{(B)} \times m}; \quad E^{(C)}: \boldsymbol{C} \to \boldsymbol{C}^{\ell_f^{(C)} \times n}; \quad E^{(D)}: \boldsymbol{C} \to \boldsymbol{C}^{\ell_f^{(D)} \times m} \quad (11)$$

The parametrical multi perturbations in the dynamic system (1)-(2) are defined by the matrices (7), subject to (8)-(11), dependent on the argument $z := \left( z^{(A)}, z^{(B)}, z^{(C)}, z^{(D)} \right)$ which takes values in some domain $\boldsymbol{C}_\alpha$ of $\boldsymbol{C}^q$. The matrices (11) are scaling matrices common to all the output components being independent of the various subscripts (i, j). The matrices in (9) are also scaling matrices of the state versus state components, state versus input component, output versus state component s and output versus input components. The matrices in (10) are specific parametrical perturbations which become weighted by the contribution of the corresponding component of z in $\boldsymbol{C}_\alpha$ through the global parametrical perturbations (7)-(8). Note by direct inspection that if all the perturbation matrices in (10) are zero then, the dynamic system (1)-(2) becomes the nominal one. Note also that an extension of the parametrical perturbations consisting of considering the scaling matrices Eqs. 11 to be dependent on the indices (i, j) would not become more general than that given in view of the whole structure of the multi perturbations (7)-(8). The following matrices are defined for each system (1)-(2) (see [1-3], [8-9]):

**Definitions**

**3.1** The spectral controllability matrix function is defined by:
$$Z_C\left(s, z^{(A)}, z^{(B)}\right) := \left(sI_n - A\left(z^{(A)}\right) - \tilde{A}\left(z^{(A)}\right) : B\left(z^{(B)}\right) + \tilde{B}\left(z^{(B)}\right)\right)$$
with that of the nominal system being $Z_{C0}\left(s, z^{(A)}, z^{(B)}\right) := \left(sI_n - A\left(z^{(A)}\right) : B\left(z^{(B)}\right)\right)$.

**3.2** The spectral observability matrix function is defined by:
$$Z_O\left(s, z^{(A)}, z^{(C)}\right) := \left(sI_n - A^*\left(z^{(A)}\right) - \tilde{A}^*\left(z^{(A)}\right) : C^*\left(z^{(C)}\right) + \tilde{C}^*\left(z^{(C)}\right)\right)^*$$
with that of the nominal system being $Z_{O0}\left(s, z^{(A)}, z^{(B)}\right) := \left(sI_n - A^*\left(z^{(A)}\right) : C^*\left(z^{(C)}\right)\right)^*$.

**3.3** The spectral output controllability matrix function is defined by:
$$Z_{OC}(s, z) := \left(C\left(z^{(C)}\right) + \tilde{C}\left(z^{(C)}\right)\right)\left(sI_n - A\left(z^{(A)}\right) - \tilde{A}\left(z^{(A)}\right) : B\left(z^{(B)}\right) + \tilde{B}\left(z^{(B)}\right)\right) + D\left(z^{(D)}\right) + \tilde{D}\left(z^{(D)}\right)$$
with that of the nominal system being $Z_{OC0}\left(s, z^{(A)}, z^{(B)}\right) := \left(sI_n - A\left(z^{(A)}\right) : B\left(z^{(B)}\right)\right)$.

**3.4** The system matrix function is defined by :

$$S(s,z) := \begin{pmatrix} sI_n - A(z^{(A)}) - \tilde{A}(z^{(A)}) & \vdots & -B(z^{(B)}) - \tilde{B}(z^{(B)}) \\ \cdots\cdots\cdots & & \cdots\cdots\cdots \\ C(z^{(C)}) + \tilde{C}(z^{(C)}) & \vdots & D(z^{(D)}) + \tilde{D}(z^{(D)}) \end{pmatrix}$$

with that of the nominal system being $S_0(s,z) := \begin{pmatrix} sI_n - A(z^{(A)}) & \vdots & -B(z^{(B)}) \\ \cdots\cdots\cdots & & \cdots\cdots \\ C(z^{(C)}) & \vdots & D(z^{(D)}) \end{pmatrix}$. □

Note by direct inspection that these matrices depend on the nominal system and the parametrical perturbations as follows

$$Z_C(s, z^{(A)}, z^{(B)}) = Z_{C0}(s, z^{(A)}, z^{(B)}) + \left(-\tilde{A}(z^{(A)}) \vdots \tilde{B}(z^{(B)})\right) \tag{12}$$

$$Z_O(s, z^{(A)}, z^{(C)}) = Z_{O0}(s, z^{(A)}, z^{(C)}) + \left(-\tilde{A}^*(z^{(A)}) \vdots \tilde{C}^*(z^{(C)})\right)^* \tag{13}$$

$$Z_{OC}(s,z) = Z_{OC0}(s,z) + C(z^{(C)})\left(-\tilde{A}(z^{(A)}) \vdots \tilde{B}(z^{(B)})\right) + \tilde{C}(z^{(C)})\left(sI_n - A(z^{(A)}) - \tilde{A}(z^{(A)}) \vdots B(z^{(B)}) + \tilde{B}(z^{(B)})\right) + \tilde{D}(z^{(D)}) \tag{14}$$

$$S(s,z) := S_0(s,z) + \begin{pmatrix} sI_n - \tilde{A}(z^{(A)}) & \vdots & -\tilde{B}(z^{(B)}) \\ \cdots\cdots\cdots & & \cdots\cdots \\ \tilde{C}(z^{(C)}) & \vdots & \tilde{D}(z^{(D)}) \end{pmatrix} \tag{15}$$

Related to Definitions 3.1-3-4 are the following ones:

**Definitions**

**3.5** $s_0 \in C$ is an input-decoupling zero of (1)-(2) for a given $(z^{(A)}, z^{(B)}) \in C^{q_A + q_B + 2}$ if $rank\ Z_C(s_0, z^{(A)}, z^{(B)}) < n$.

**3.6** $s_0 \in C$ is an output-decoupling zero of (1)-(2) for a given $(z^{(A)}, z^{(C)}) \in C^{q_A + q_C + 2}$ if $rank\ Z_O(s_0, z^{(A)}, z^{(C)}) < n$.

**3.7** $s_0 \in C$ is an input/output-decoupling zero of (1)-(2) for a given $z \in C^q$ if $max\left(rank\ Z_C(s_0, z^{(A)}, z^{(B)}), rank\ Z_O(s_0, z^{(A)}, z^{(C)})\right) < n$.

**3.8** $s_0 \in C$ is an external input-decoupling zero of (1)-(2) for a given $z \in C^q$ if $rank\ Z_{OC}(s_0, z) < p$.

**3.9** $s_0 \in C$ is an invariant zero of (1)-(2) for a given $z \in C^q$ if $rank\ S(s_0, z) < n + min(m, p)$.

**3.10** $s_0 \in C$ is a transmission zero of (1)-(2) for a given $z \in C^q$ if it is an invariant zero which is not an input-decoupling or output-decoupling zero; i.e. $rank\ S(s_0, z) < n + min(m, p)$ and $rank\ Z_C(s_0, z^{(A)}, z^{(B)}) = rank\ Z_O(s_0, z^{(A)}, z^{(C)}) = n$. □

It is well-known that the system (1)-(2) is controllable (respectively, observable) for a certain $z \in C^q$ if it has no input-decoupling zero (respectively, no output-decoupling zero). Invariant zeros which are not decoupling zeros are transmission zeros in the sense that if $s_0 \in C$ is a transmission zero for a certain $z \in C^q$ then $y \equiv 0$ if $u(t) = K_u e^{s_0 t}$ for any $K_u \in R$ if $x_0 = 0$ (input-output transmission blocking property) and

$$\hat{G}(s,z) = \left(C(z) + \tilde{C}(z)\right)\left(sI_n - A(z^{(A)}) - \tilde{A}(z^{(A)})\right)^{-1}\left(B(z^{(B)}) + \tilde{B}^{(B)}\right) + \left(D(z^{(D)}) + \tilde{D}(z^{(D)})\right)$$

is the transfer matrix of the system (1)-(2) defined as $\hat{G}(s,z) := \hat{y}(s,z)/\hat{u}(s,z)$ for $x_0 = 0$. Note that decoupling zeros are poles of the system transfer matrix. The transmission zeros are zeros of $\hat{G}(s,z)$ which are not poles of (1)-(2), i.e. which are not eigenvalues of $A(z) + \tilde{A}(z)$. Note also that if the system is controllable and observable all invariant zero, if any, is a transmission zero. Input/output decoupling zeros are poles cancelled by zeros in the transfer function so that they are not transmission zeros. Finally, note that input/output-decoupling zeros are invariant zeros since

$$max\left(rank\, Z_C\left(s_0, z^{(A)}, z^{(B)}\right), rank\, Z_O\left(s_0, z^{(A)}, z^{(C)}\right)\right) < n \Rightarrow rank\, S(s_0, z) < n + min(m, p)$$

However, input-decoupling zeros (respectively, output-decoupling zeros) which are not output-decoupling zeros (respectively, input-decoupling zeros) are not invariant zeros since

$$rank\, Z_C\left(s_0, z^{(A)}, z^{(B)}\right) < n \text{ may imply } rank\, S(s_0, z) = n + m \text{ for } p > m$$
$$rank\, Z_O\left(s_0, z^{(A)}, z^{(C)}\right) < n \text{ may imply } rank\, S(s_0, z) = n + p \text{ for } m > p$$

The system (1)-(2) is said to be controllable if there is a control vector function on $[0, T]$ such that the state takes any prescribed finite value at any finite prescribed time T. It is said to be observable if any bounded initial can be computed from measures of the output vector on any finite time interval $[0, T]$. Those properties coincide in the linear-time-invariant case with the respective properties of spectral controllability/spectral observability which hold if and only if the spectral controllability/spectral observability matrix functions are full rank $\forall (s, z) \in C \times C_\alpha$. Thus, spectral controllability/ observability properties will be referred to in the following simply as controllability/observability. Definitions 3.5 to 3.10 combined with Popov-Belevitch-Hautus controllability and observability tests [1] lead to the subsequent result which considers bounded sets where the varying parameters belong to defined by $C_\alpha^{(A)} \subset \{1\} \times C^{q_A}$, $C_\alpha^{(B)} \subset \{1\} \times C^{q_B}$, $C_\alpha^{(C)} \subset \{1\} \times C^{q_C}$ and $C_\alpha^{(D)} \subset \{1\} \times C^{q_D}$ :

**Theorem 3.11.** The following properties hold**:**

**(i)** The system (1)-(2) is controllable in a bounded domain $C_\alpha^{(A,B)} := C_\alpha^{(A)} \times C_\alpha^{(B)} \subset \mathbf{C}^{q_A + q_B + 2}$, if and only if it has no input-decoupling zero in $C_\alpha^{(A,B)}$, and equivalently, if and only if

$$\text{rank } Z_C\left(s, z^{(A)}, z^{(B)}\right) = n, \forall \left(s, z^{(A)}, z^{(B)}\right) \in Im(s_f) \times C_\alpha^{(A,B)} \subset \mathbf{C}^{q_A + q_B + 3}$$

where the discrete function $s_f$ is defined as $s_f : C_\alpha^{(A)} \to \sigma\left(A(z^{(A)}) + \tilde{A}(z^{(A)})\right)$.

**(ii)** The system (1)-(2) is observable in a bounded domain $C_\alpha^{(A,C)} := C_\alpha^{(A)} \times C_\alpha^{(C)} \subset \mathbf{C}^{q_A + q_C + 2}$ if and only if it has no output-decoupling zero in $C_\alpha^{(A,C)}$, and equivalently, if and only if

$$\text{rank } Z_O\left(s, z^{(A)}, z^{(C)}\right) = n, \forall \left(s, z^{(A)}, z^{(C)}\right) \in Im(s_f) \times C_\alpha^{(A,C)} \subset \mathbf{C}^{q_A + q_C + 3}.$$

**(iii)** The system (1)-(2) is controllable and observable in a bounded domain $C_\alpha^{(A,B,C)} := C_\alpha^{(A)} \times C_\alpha^{(B)} \times C_\alpha^{(C)} \subset \mathbf{C}^{q_A + q_B + q_C + 3}$ if and only if it has no input-decoupling zero in $C_\alpha^{(A)} \times C_\alpha^{(B)}$ and no output-decoupling zero in $C_\alpha^{(A)} \times C_\alpha^{(C)}$, and equivalently, if and only if

$$\text{rank } Z_C\left(s, z^{(A)}, z^{(B)}\right) = \text{rank } Z_O\left(s, z^{(A)}, z^{(C)}\right) = n,$$

$\forall \left(s, z^{(A)}, z^{(B)}, z^{(C)}\right) \in Im(s_f) \times C_\alpha^{(A,B,C)} \subset \mathbf{C}^{q_A + q_B + q_C + 3}$. As a result, The system (1)-(2) is controllable and observable in $C_\alpha^{(A,B,C)} \Leftrightarrow \neg \exists$ a system (1)-(2) with $dim(x(t)) < n$ for some $z \in C_\alpha$ and transfer matrix

$$\hat{G}(s,z) = \left(C(z) + \tilde{C}(z)\right)\left(sI_n - A(z^{(A)}) - \tilde{A}(z^{(A)})\right)^{-1}\left(B(z^{(B)}) + \tilde{B}^{(B)}\right) + \left(D(z^{(D)}) + \tilde{D}(z^{(D)})\right)$$

**(iv)** The system (1)-(2) is output controllable in a bounded domain $C_\alpha := C_\alpha^{(A)} \times C_\alpha^{(B)} \times C_\alpha^{(C)} \times C_\alpha^{(D)}$ if and only if it has no external input-decoupling zero in $C_\alpha$, and equivalently, if and only if

$$\text{rank } Z_{OC}\left(s, z^{(A)}, z^{(B)}, z^{(C)}, z^{(D)}\right) = p, \forall \left(s, z^{(A)}, z^{(B)}, z^{(C)}, z^{(D)}\right) \in \mathbf{C} \times C_\alpha \subset \mathbf{C}^{q_A + q_B + q_C + q_D + 5}.$$

If $rank\left(C(z) + \tilde{C}(z)\right) = p$ in $C_\alpha^{(A)} \times C_\alpha^{(B)} \times C_\alpha^{(C)}$ then the system (1)-(2) is output controllable if and only if

$$\text{rank } Z_{OC}\left(s, z^{(A)}, z^{(B)}, z^{(C)}\right) = p$$

$, \forall \left(s, z^{(A)}, z^{(B)}, z^{(C)}\right) \in Im(s_f) \times C_\alpha^{(A)} \times C_\alpha^{(B)} \times C_\alpha^{(C)} \subset \mathbf{C}^{q_A + q_B + q_C + 3}$

**Proof: (i)** The continuous vector function $s_f : C_\alpha^{(A)} \to \sigma\left(A(z^{(A)}) + \tilde{A}(z^{(A)})\right)$ exists on $C_\alpha^{(A)}$ since the eigenvalues of the square matrix function $A(z^{(A)}) + \tilde{A}(z^{(A)})$ exist and are bounded continuous functions on the bounded domain $C_\alpha^{(A)}$. Define the graph of $s_f$ on its definition domain as $G(s_f) := \left\{ (\omega, s_f(\omega)) : \omega \in C_\alpha^{(A)} \right\}$. Taking Laplace transforms in (1) gives the equivalent algebraic linear equation

$$Z_C(s, z^{(A)}, z^{(B)})\left(\hat{x}^*(s, z^{(A)}, z^{(B)}) \vdots -\hat{u}^*(s, z^{(A)}, z^{(B)})\right)^* = x_0$$

for any bounded initial conditions $x(0, z^{(A)}, z^{(B)}) = x_0$. Take a linear state-feedback control of Laplace transform $\hat{u}(s, z^{(A)}, z^{(B)}) = K(s, z^{(A)}, z^{(B)}) \hat{x}(s, z^{(A)}, z^{(B)})$ for some $K : C \times C_\alpha^{(A)} \times C_\alpha^{(B)} \subset C \times C^{q_A} \times C^{q_B} \to range(K) \subset C^{m \times n}$ is a matrix function of order $m \times n$. Combining the above equations, one gets the linear algebraic system:

$$Z_C(s, z^{(A)}, z^{(B)})\left(I_n \vdots -K^*(s, z^{(A)}, z^{(B)})\right)^* \hat{x}(s, z^{(A)}, z^{(B)}) = x(0, z^{(A)}, z^{(B)})$$

Since $rank\, Z_C(s, z^{(A)}, z^{(B)}) = n$ then the generic rank of the square n-matrix function $Z_C(s, z^{(A)}, z^{(B)})\left(I_n \vdots -K^*(s, z^{(A)}, z^{(B)})\right)^*$, considered as a function of the feedback matrix $K(s, z^{(A)}, z^{(B)})$, is n everywhere in $C \times C_\alpha^{(A,B)}$ so that the above algebraic linear system may be full rank at any point of the domain of $K(s, z^{(A)}, z^{(B)})$ if its range is (pointwise) chosen appropriately. As a result,

$$rank\left[Z_C(s, z^{(A)}, z^{(B)})\left(I_n \vdots -K^*(s, z^{(A)}, z^{(B)})\right)^*\right] < n$$
$$\Leftrightarrow det\left(sI_n - A(z^{(A)}) - \tilde{A}(z^{(A)}) - (B(z^{(B)}) + \tilde{B}(z^{(B)}))K(s, z^{(A)}, z^{(B)})\right) = 0$$

at arbitrarily fixed complex points $s(z^{(A)}, z^{(B)})$ in $C_\alpha^{(A,B)}$ since any $s \in C$ may be chosen not to be an eigenvalue of the matrix

$$A(z^{(A)}) + \tilde{A}(z^{(A)}) + (B(z^{(B)}) + \tilde{B}(z^{(B)}))K(s, z^{(A)}, z^{(B)})$$

everywhere in $C_\alpha^{(A,B)}$. Since the eigenvalues are arbitrarily assignable by linear state-feedback the system (1)-(2) is controllable if $rank\, Z_C(s, z^{(A)}, z^{(B)}) = n$ what proves sufficiency. To prove necessity, proceed by contradiction by assuming that $rank\, Z_C\left(s, z^{(A)}, z^{(B)}\right) < n$. Then, there exists a vector $0 \ne q \in \boldsymbol{C}^n$ for each point in $C_\alpha^{(A,B)}$ such that :

$$q^* Z_C\left(s, z^{(A)}, z^{(B)}\right) = 0 \Leftrightarrow q^*\left(sI_n - A\left(z^{(A)}\right) - \tilde{A}\left(z^{(A)}\right)\right) = q^*\left(B\left(z^{(B)}\right) + \tilde{B}\left(z^{(B)}\right)\right) = 0$$

i.e. q is a nonzero eigenvector of $A\left(z^{(A)}\right) + \tilde{A}\left(z^{(A)}\right)$, with associate eigenvalue $\lambda_q$, which is orthogonal to the matrix $B\left(z^{(B)}\right) + \tilde{B}\left(z^{(B)}\right)$. Direct calculation with (1)-(2) yields that a state trajectory solution satisfies $\dot{x}\left(t, z^{(A)}, z^{(B)}\right) = \lambda_q x\left(t, z^{(A)}, z^{(B)}\right)$ with associate state trajectory solution $x\left(t, z^{(A)}, z^{(B)}\right) = e^{\lambda_q t} q$ irrespective of the control, $\forall t \in \boldsymbol{R}_{+0}$. It is obvious that such a trajectory cannot reach any point $\bar{x} = x\left(T, z^{(A)}, z^{(B)}\right) \ne e^{\lambda_q T} q$ at any finite arbitrary time $T > 0$. Thus, the system is not controllable in $C_\alpha^{(A,B)}$ which proves necessity. Thus, the system (1)-(2) is controllable for any $z \in C_\alpha^{(A,B)}$ if and only if

$$rank\, Z_C\left(s, z^{(A)}, z^{(B)}\right) = n,\, \forall\left(s, z^{(A)}, z^{(B)}\right) \in \boldsymbol{C} \times C_\alpha^{(A,B)}$$
$$\Leftrightarrow rank\, Z_C\left(s, z^{(A)}, z^{(B)}\right) = n,\, \forall\left(s, z^{(A)}, z^{(B)}\right) \in G\left(s_f\right) \times C_\alpha^{(B)}$$
$$\Leftrightarrow rank\, Z_C\left(\omega, z^{(A)}, z^{(B)}\right) = n,\, \forall\left(s, z^{(A)}, z^{(B)}\right) \in Im\left(s_f\right) \times C_\alpha^{(A,B)}$$

since loss of rank in some point of $\boldsymbol{C} \times C_\alpha^{(A,B)}$ is only possible for $s \in \sigma\left(A\left(z^{(A)}\right) + \tilde{A}\left(z^{(A)}\right)\right)$, which is $Im\left(s_f\right)$, from Popov-Belevitch-Hautus rank controllability test. Equivalently, the system (1)-(2) is controllable if and only if it has no input-decoupling zero in $C_\alpha^{(A,B)}$ from Definition 3.5.

**(ii)** The proof is similar from the Popov-Belevitch-Hautus observability rank test $rank\, Z_O\left(s, z^{(A)}, z^{(C)}\right) = n,\, \forall\left(s, z^{(A)}, z^{(C)}\right) \in G\left(s_f\right) \times C_\alpha^{(C)}$, and, equivalently, $\forall\left(s, z^{(A)}, z^{(C)}\right) \in Im\left(s_f\right) \times C_\alpha^{(A,C)}$, since observability is a dual property to controllability through the replacements $A \to A^*$, $B \to C^*$.

**(iii)** The first part follows by combining the proofs of Properties (i)-(ii). The second part is now proven. Assume that there is at least a transmission zero $s_0 \in \boldsymbol{C}$ and $S(s_0, z)\left(\hat{x}^T(s), u^T(s)\right)^T = 0 \Rightarrow rank\, S(s_0, z) < n$ for some $\left(\hat{x}(s_0, z), \hat{u}^T(s_0, z)\right)^T \ne (0, 0)^T \in \boldsymbol{R}^{n+m}$ and some $z \in C_\alpha$ since $s_0 \in \boldsymbol{C}$ is neither an input-decoupling or and output decoupling zero of (1)-(2). Since $s_0 \notin \sigma\left(A(z) + \tilde{A}(z)\right)$

$$\det S(s_0, z) = \det\left(\operatorname{diag}\left(s_0 I_n - A(z^{(A)}) - \tilde{A}(z^{(A)}) \;\vdots\; I_n\right)\right)$$

$$\times \det\begin{pmatrix} I_n & \vdots & -\left(s_0 I_n - A(z^{(A)}) - \tilde{A}(z^{(A)})\right)^{-1} B(z^{(B)}) + \tilde{B}(z^{(B)}) \\ \cdots & \cdots & \cdots \\ C(z^{(C)}) + \tilde{C}(z^{(C)}) & \vdots & D(z^{(D)}) + \tilde{D}(z^{(D)}) \end{pmatrix} = 0$$

$\Leftrightarrow C(z^{(C)}) + \tilde{C}(z^{(C)}) \operatorname{Adj}\left(s_0 I_n - A(z^{(A)}) - \tilde{A}(z^{(A)})\right) B(z^{(B)}) + \tilde{B}(z^{(B)}) + D(z^{(D)}) + \tilde{D}(z^{(D)}) = 0$

$\Leftrightarrow \hat{y}(s_0, z) \equiv 0$ for some $\hat{u}(s_0, z) \neq 0 \in \mathbf{R}^m$ and some $z \in C_\alpha$ if $x_0 = 0$ and the system possesses the input-output transmission blocking property at the transmission zero $s = s_0$ for some $z \in C_\alpha^{(A,B,C)}$. Since $s_0 \notin \sigma\left(A(z) + \tilde{A}(z)\right)$, $\hat{G}(s, z)$ has no zero-pole cancellation at $s = s_0$ so that there is no system (1)-(2) with state dimension less than n which possess the input-output transmission blocking property for some $z \in C_\alpha^{(A,B,C)}$.

**(iv)** Its proof is similar to those of (i)-(ii) from the Popov-Belevitch-Hautus output controllability rank test, namely, $\operatorname{rank} Z_{OC}\left(s, z^{(A)}, z^{(B)}, z^{(C)}, z^{(D)}\right) = p$, $\forall s \in \mathbf{C} \times C_\alpha$. $\square$

Note that Theorem 3.11 also holds if $C_\alpha$ is a closed domain. In this case, $G(s_f)$ is closed since there is a finite number of eigenvalues of each matrix $A(z^{(A)}) + \tilde{A}(z^{(A)})$ in $C_\alpha^{(A)}$. As a result, all pairs in the set $G(s_f)$ conform a closed set since $C_\alpha^{(A)}$ is closed. The subsequent results is a direct consequence of Theorem 3.11. It is of interest to formulate the various properties for the system (1)-(2) under parametrical multi perturbations belonging to a certain domain in an easy testable form provided that provided that they hold for the nominal system in an easy testable form. Note also that Theorem 3.11 is not directly applicable to time-varying parameters but to varying parameterizations within some appropriate domains.

**Corollary 3.12.** The following properties hold:

(i) The system (1)-(2) is controllable in $C_\alpha^{(A,B)}$ if and only if

$$\det \hat{Z}_C\left(s, z^{(A)}, z^{(B)}\right) > 0, \forall \left(s, z^{(A)}, z^{(B)}\right) \in \mathbf{C} \times C_\alpha^{(A,B)}$$

$$\Leftrightarrow \operatorname{Inf}\left(\underline{\sigma}\left(Z_C\left(s, z^{(A)}, z^{(B)}\right)\right): \left(s, z^{(A)}, z^{(B)}\right) \in \mathbf{C} \times C_\alpha^{(A,B)}\right) > 0$$

where $\hat{Z}_C\left(s, z^{(A)}, z^{(B)}\right) := Z_C\left(s, z^{(A)}, z^{(B)}\right) Z_C^*\left(s, z^{(A)}, z^{(B)}\right)$ and, equivalently, if and only if

$$\det \hat{Z}_C\left(s, z^{(B)}\right) > 0, \forall \left(s, z^{(B)}\right) \in G(s_f) \times C_\alpha^{(B)}$$

$$\Leftrightarrow Inf\left(\underline{\sigma}\left(Z_C\left(s, z^{(A)}, z^{(B)}\right)\right):\left(s, z^{(B)}\right) \in G\left(s_f\right) \times C_\alpha^{(B)}\right) > 0.$$

**(ii)** The system (1)-(2) is observable in $C_\alpha$ if and only if

$$det\, \hat{Z}_O\left(s, z^{(A)}, z^{(C)}\right) > 0, \forall\left(s, z^{(A)}, z^{(C)}\right) \in C \times C_\alpha^{(A,C)}$$

$$\Leftrightarrow Inf\left(\underline{\sigma}\left(Z_O\left(s, z^{(A)}, z^{(C)}\right)\right):\left(s, z^{(A)}, z^{(C)}\right) \in C \times C_\alpha^{(A,C)}\right) > 0$$

where $\hat{Z}_O\left(s, z^{(A)}, z^{(C)}\right) := Z_O^*\left(s, z^{(A)}, z^{(C)}\right) Z_O\left(s, z^{(A)}, z^{(C)}\right)$ and, equivalently, if and only if

$$det\, \hat{Z}_O\left(s, z^{(A)}, z^{(C)}\right) > 0, \forall\left(s, z^{(C)}\right) \in G\left(s_f\right) \times C_\alpha^{(C)}$$

$$\Leftrightarrow Inf\left(\underline{\sigma}\left(Z_O\left(s, z^{(A)}, z^{(C)}\right)\right):\left(s, z^{(C)}\right) \in G\left(s_f\right) \times C_\alpha^{(C)}\right) > 0.$$

**(iii)** The system (1)-(2) is controllable and observable in $C_\alpha^{(A,B,C)}$ if and only if

$$det\, \hat{Z}_C\left(s, z^{(A)}, z^{(B)}\right) \cdot det\, \hat{Z}_O\left(s, z^{(A)}, z^{(C)}\right) > 0 \Leftrightarrow \underline{\sigma}\left(Z_C\left(s, z^{(A)}, z^{(B)}\right)\right) \cdot \underline{\sigma}\left(Z_O\left(s, z^{(A)}, z^{(C)}\right)\right) > 0$$

$$, \forall\left(s, z^{(B)}, z^{(C)}\right) \in G\left(s_f\right) \times C_\alpha^{(B,C)}.$$

**(iv)** The system (1)-(2) is output controllable in $C_\alpha$ if and only if $Z_{OC}(s, z)$ has no zero singular value $\forall z \in C_\alpha$ and, equivalently, if and only if $det\, \hat{Z}_{OC}(s, z) > 0$, $\forall(s,z) \in C \times C_\alpha$, where $\hat{Z}_{OC}(s, z) := Z_{OC}(s, z) Z_{OC}^*(s, z)$.

**Proof**: **(i)** Note that $\hat{Z}_C\left(s, z^{(A)}, z^{(B)}\right)$ is a (n+m)- square matrix which is normal by construction. Thus, its eigenvalues are nonnegative and real being the squares of the singular values of $Z_C\left(s, z^{(A)}, z^{(B)}\right)$. Consider again the function $s_f : C_\alpha^{(A)} \to \sigma\left(A\left(z^{(A)}\right) + \tilde{A}\left(z^{(A)}\right)\right)$ of Theorem 3.11. Thus, the system (1)-(2) is controllable from Theorem 3.11(i) if and only if

$$rank\, Z_C\left(s, z^{(A)}, z^{(B)}\right) = n, \forall\left(s, z^{(A)}, z^{(B)}\right) \in C \times C_\alpha^{(A,B)}$$

$$\Leftrightarrow rank\, Z_C\left(s, z^{(B)}\right) = n, \forall\left(s, z^{(B)}\right) \in G\left(s_f\right) \times C_\alpha^{(B)}$$

$$\Leftrightarrow rank\, \hat{Z}_C\left(s, z^{(A)}, z^{(B)}\right) = n, \forall\left(s, z^{(A)}, z^{(B)}\right) \in C \times C_\alpha^{(A,B)}$$

$$\Leftrightarrow det\, \hat{Z}_C\left(s, z^{(A)}, z^{(B)}\right) > 0, \forall\left(s, z^{(A)}, z^{(B)}\right) \in C \times C_\alpha^{(A,B)}$$

$$\Leftrightarrow rank\, \hat{Z}_C\left(\omega, z^{(B)}\right) = n, \forall\left(\omega, z^{(B)}\right) \in G\left(s_f\right) \times C_\alpha^{(B)}$$

$$\Leftrightarrow det\, \hat{Z}_C\left(\omega, z^{(B)}\right) > 0, \forall\left(\omega, z^{(B)}\right) \in G\left(s_f\right) \times C_\alpha^{(B)}$$

$$\Leftrightarrow \underline{\sigma}\left(Z_C\left(s, z^{(A)}, z^{(B)}\right)\right) > 0, \forall\left(s, z^{(A)}, z^{(B)}\right) \in C \times C_\alpha^{(A,B)}$$

$$\Leftrightarrow \underline{\sigma}\left(Z_C\left(\omega, z^{(B)}\right)\right) > 0, \forall \left(\omega, z^{(B)}\right) \in G\left(s_f\right) \times C_\alpha^{(B)}$$

what proves Property (i). The proofs of (ii)-(iv) are similar from the parallel properties of Theorem 3.11 and are then omitted. □

**Remark 3.13.** Since controllability is lost in $C_\alpha$ if and only if $Z_C\left(s, z^{(B)}\right)$ is rank defective for some $\left(s, z^{(B)}\right) \in G\left(s_f\right) \times C_\alpha^{(B)}$ [Theorem 3.11 (i)], it turns out that controllability in $C_\alpha$ holds if and only if $\underline{\sigma}\left(Z_C\left(s, z^{(B)}\right)\right) > 0, \forall \left(s, z^{(B)}\right) \in G\left(s_f\right) \times C_\alpha^{(B)}$ and, equivalently, if and only if $\det \hat{Z}_C\left(s, z^{(B)}\right) > 0, \forall \left(s, z^{(B)}\right) \in G\left(s_f\right) \times C_\alpha^{(B)}$. As a result, the tests of Corollary 3.12 (i) have only to be performed for $s \in \sigma\left(A(z) + \tilde{A}(z)\right)$ for each $\left(z^{(A)}, z^{(B)}\right) \in C_\alpha^{(A,B)}$. Similar considerations apply to Corollary 3.12 [(ii)-(iv)]. □

Theorem 3.11 and Corollary 3.12 are directly extendable to stabilizability and detectability as follows. First, define $s_f^{+0}: Dom\left(s_f^{+0}\right) \to \sigma\left(A\left(z^{(A)}\right) + \tilde{A}\left(z^{(A)}\right)\right) \cap \boldsymbol{C}_{+0}$ of graph

$$G\left(s_f^{+0}\right) := \left\{\left(\omega, s_f^{+0}(\omega)\right) : \omega \in Dom\left(s_f^{+0}\right)\right\}$$ where

$$C_\alpha^{(A)} \supset Dom\left(s_f^{+0}\right) := \left\{\omega \in C_\alpha^{(A)} : s_f^{+0}(\omega) \in \boldsymbol{C}_{+0}\right\}$$

Note that $G\left(s_f^{+0}\right) \subset G\left(s_f\right)$ since $Dom\left(s_f^{+0}\right) \subset Dom\left(s_f\right) \equiv C_\alpha^{(A)}$ and there is a natural projection function $p_f^{+0}: Dom\left(p_f^{0+}\right) \equiv Im(s_f) \to Im\left(s_f^{+0}\right) \subset Im(s_f)$.

**Corollary 3. 14**.

**(i)** The system (1)-(2) is stabilizable in a bounded domain $C_\alpha$ if and only any of the following equivalent properties hold:

**(i.1)** All the input-decoupling zeros in $C_\alpha^{(A,B)}$, if any, have negative real parts

**(i.2)** $rank\, Z_C\left(s, z^{(A)}, z^{(B)}\right) = n, \forall \left(s, z^{(A)}, z^{(B)}\right) \in Im\left(s_f^{+0}\right) \times C_\alpha^{(A,B)}$

**(i.3)** $rank\, Z_C\left(s, z^{(A)}, z^{(B)}\right) = n, \forall \left(s, z^{(A)}, z^{(B)}\right) \in \boldsymbol{C}_{+0} \times C_\alpha^{(A,B)}$

**(i.4)** $\det \hat{Z}_C\left(s, z^{(B)}\right) > 0, \forall \left(s, z^{(B)}\right) \in G\left(s_f^{+0}\right) \times C_\alpha^{(B)}$

**(i.5)** $\underline{\sigma}\left(Z_C\left(G\left(s, z^{(B)}\right)\right)\right) > 0, \forall \left(s, z^{(B)}\right) \in G\left(s_f^{+0}\right) \times C_\alpha^{(B)}$

**(ii)** The system (1)-(2) is detectable in $C_\alpha$ if and only any of the following equivalent

properties hold:

**(ii.1)** All the output-decoupling zeros in $C_\alpha^{(A,C)}$, if any, have negative real parts

**(ii.2)** $rank\ Z_O\left(s, z^{(A)}, z^{(C)}\right) = n, \forall \left(s, z^{(A)}, z^{(C)}\right) \in C_{+0} \times C_\alpha^{(A,C)}$.

**(ii.3)** $rank\ Z_O\left(s, z^{(A)}, z^{(C)}\right) = n, \forall \left(s, z^{(A)}, z^{(C)}\right) \in Im\left(s_f^{+0}\right) \times C_\alpha^{(A,C)}$.

**(ii.4)** $det\ \hat{Z}_O\left(s, z^{(C)}\right) > 0, \forall \left(s, z^{(C)}\right) \in G\left(s_f^{+0}\right) \times C_\alpha^{(C)}$

**(ii.5)** $\underline{\sigma}\left(Z_C\left(G\left(s, z^{(C)}\right)\right)\right) > 0, \forall \left(s, z^{(C)}\right) \in G\left(s_f^{+0}\right) \times C_\alpha^{(C)}$ □

**Remarks 3.15.** The system (1)-(2) is said to be stabilizable if there is no input-decoupling zero in $C_{0+}$. If the system (1)-(2) is stabilizable then there exists some state-feedback control law $u: R_{0+} \times X \to U$ such that the system (1)-(2) is globally asymptotically Lyapunov´s stable under such a law. The system (1)-(2) is said to be detectable if there is no output-decoupling zero in $C_{+0}$. From Theorem 3.12 and Corollaries 3.12 and 3.14, it turns out that if the system (1)-(2) is controllable (respectively, observable) in $C_\alpha$ then it is stabilizable (respectively, detectable) in $C_\alpha$. Also, if the unforced system (1)-(2) is globally asymptotically stable in the Lyapunov´s sense on $C_\alpha$ (i.e. the matrix $A\left(z^{(A)}\right) + \tilde{A}\left(z^{(A)}\right)$ has all its eigenvalues with negative real parts for all $z^{(A)} \in C_\alpha^{(A)}$ then it is also stabilizable (even if it is not controllable) and detectable (even if it is not observable) since the spectral controllability, respectively, observability matrices are jointly full rank on $C_{+0} \times C_\alpha^{(A,B)}$, respectively, $C_{+0} \times C_\alpha^{(A,B)}$. □

**Remarks 3.16.** Note that the conditions implying the corresponding determinants or minimum singular values to be positive to guarantee each of the controllability, stabilizability, observability and detectability properties in Theorem 3.11 and Corollaries 3.13 and 3.14 is sufficient in a bounded domain. However, the boundedness of determinants and maximum singular values is also needed to guarantee each property in an unbounded domain.

## 4. Maintenance of the properties from those of the nominal system

In this section, attention is paid to the normal matrix $\hat{Z}_C\left(s, z^{(A)}, z^{(B)}\right) = Z_C\left(s, z^{(A)}, z^{(B)}\right) Z_C^*\left(s, z^{(A)}, z^{(B)}\right)$ to obtain conditions for maintaining or loosing controllability under parametrical multi perturbations of a certain size provided that the nominal system is controllable. In the analysis, it is taken advantage of the fact that such a matrix is square and nonsingular if the system is controllable. Furthermore, if the matrix $\left(A\left(z^{(A)}\right) + \tilde{A}\left(z^{(A)}\right): B\left(z^{(B)}\right) + \tilde{B}\left(z^{(B)}\right)\right)$ is Hermitian, so that $\left(sI_n - \left(A\left(z^{(A)}\right) + \tilde{A}\left(z^{(A)}\right)\right): B\left(z^{(B)}\right) + \tilde{B}\left(z^{(B)}\right)\right)$ is also Hermitian, then the eventual loss of rank of $\hat{Z}_C(s, z^{(A)}, z^{(B)})$ for some $(z^{(A)}, z^{(B)})$ only occurs for $s \in \mathbb{C}$ being some singular value of $Z_C(s, z^{(A)}, z^{(B)})$, i.e. for some eigenvalue of $A(z^{(A)}) + \tilde{A}(z^{(A)})$. In the general case, the loss of rank can occur also for eigenvalues of $A^*(z^{(A)}) + \tilde{A}^*(z^{(A)})$. A close discussion is directly applicable to stabilizability and direct extensions, by modifying accordingly the matrices, are also applicable to observability, detectability and output controllability. Using (12),

$$\hat{Z}_C(s, z^{(A)}, z^{(B)}) = \hat{Z}_{C0}(s, z^{(A)}, z^{(B)})$$
$$\times \left( I_n + \left( Z_{C0}(s, z^{(A)}, z^{(B)}) Z_{C0}^*(s, z^{(A)}, z^{(B)}) \right)^{-1} \left( \left( -\tilde{A}(z^{(A)}) \vdots \tilde{B}(z^{(B)}) \right) \left( -\tilde{A}(z^{(A)}) \vdots \tilde{B}(z^{(B)}) \right)^* \right. \right.$$
$$\left. \left. + Z_{C0}(s, z^{(A)}, z^{(B)}) \left( -\tilde{A}(z^{(A)}) \vdots \tilde{B}(z^{(B)}) \right)^* + \left( -\tilde{A}(z^{(A)}) \vdots \tilde{B}(z^{(B)}) \right) Z_{C0}^*(s, z^{(A)}, z^{(B)}) \right) \right) \quad (16)$$

provided that $\hat{Z}_{C0}(s, z^{(A)}, z^{(B)}) := Z_{C0}(s, z^{(A)}, z^{(B)}) . Z_{C0}^*(s, z^{(A)}, z^{(B)})$ is nonsingular, i.e. its minimum singular value is positive, that is, there is no zero-input decoupling zero of $Z_{C0}(s, z^{(A)}, z^{(B)})$. Expanding the first identity of (7), one gets:

$$\left( -\tilde{A}(z^{(A)}) \vdots \tilde{B}(z^{(B)}) \right) = \left( -z^{(A)} I_n \vdots z^{(B)} I_n \right) \Delta^{(AB)} \quad (17)$$

with

$$\Delta^{(AB)} := diag \left( \sum_{j=1}^{n_A} diag \left( D_{0j}^{(A)} \vdots D_{1j}^{(A)} \vdots \ldots \vdots D_{q_A j}^{(A)} \right) \left( \Delta_{0j}^{(A)*} \vdots \Delta_{1j}^{(A)*} \vdots \ldots \vdots \Delta_{q_A j}^{(A)*} \right)^* E^{(A)} \right.$$
$$\left. \vdots \sum_{j=1}^{n_B} diag \left( D_{0j}^{(B)} \vdots D_{1j}^{(B)} \vdots \ldots \vdots D_{q_B j}^{(B)} \right) \left( \Delta_{0j}^{(B)*} \vdots \Delta_{1j}^{(B)*} \vdots \ldots \vdots \Delta_{q_B j}^{(B)*} \right)^* E^{(B)} \right) \quad (18)$$

Eq. 16 may be rewritten under (17)-(18) as

$$\hat{Z}_C(s, z^{(A)}, z^{(B)}) = \hat{Z}_{C0}(s, z^{(A)}, z^{(B)}) \left( I_n + \hat{Z}_{C0}(s, z^{(A)}, z^{(B)})^{-1} \right.$$
$$\times \left( Z_{C0}(s, z^{(A)}, z^{(B)}) \Delta^{(AB)*} \left( -z^{(A)} I_n \vdots z^{(B)} I_n \right)^* + \Delta^{(AB)} Z_{C0}^*(s, z^{(A)}, z^{(B)}) \right.$$
$$\left. \left. + \left( -z^{(A)} I_n \vdots z^{(B)} I_n \right) \Delta^{(AB)} \Delta^{(AB)*} \left( -z^{(A)} I_n \vdots z^{(B)} I_n \right)^* \right) \right) \quad (19)$$

A direct result follows:

**Theorem 4.1**. The following properties hold:

**(i)** Assume that the nominal system (1)-(2) is controllable in a bounded domain $C_{C0} \subset C^{q_A+q_B+2}$. Then, there exist parametrical multi perturbations $\left(-\tilde{A}(z^{(A)}): \tilde{B}(z^{(B)})\right)$ satisfying (17)-(18), defined by maps $C_C \subset C_{C0} \to P_A \times P_B \subset C^{n \times (n+m)}$ such that the perturbed systems (1)-(2) are also controllable, and equivalently they have no input-decoupling zeros in $C$, within each given bounded domain $C_C \subset C_{C0}$. All those multi perturbations are subject to a computable norm upper-bound depending on $C_C$.

**(ii)** Assume that there is some nonzero vector $v_Q \in C^{n^2}$ which is some linear combination of the columns of the matrix

$$Q\left(z^{(A)}, z^{(B)}\right) := \left(-\sum_{i=0}^{q_A}\sum_{j=1}^{n_A} z_i^{(A)} D_{ij}^{(A)} \otimes E^{(A)^*} \ \vdots \ \sum_{i=0}^{q_B}\sum_{j=1}^{n_B} z_i^{(B)} D_{ij}^{(B)} \otimes E^{(B)^*}\right)$$

for some given $C_C \subset C_{C0}$. Then, there is a parametrical multi perturbation within the class (7) which violates the norm upper-bound referred to in Property (i) such that the system is uncontrollable in $C_C$.

**(iii)** Assume that the nominal system (1)-(2) is controllable in $C_{C0} \equiv C^{q_A+q_B+2}$ and that the scaling matrices $E^{(A)}$ and $E^{(B)}$ are both full rank. Then, there exist parametrical multi perturbations $\left(-\tilde{A}(z^{(A)}): \tilde{B}(z^{(B)})\right)$ satisfying (17)-(18), defined by maps $C^{q_A+q_B+2} \to P_A \times P_B \subset C^{n \times (n+m)}$ such that the perturbed systems (1)-(2) is not controllable in $C^{q_A+q_B+2}$.

**(iv)** Properties (i) and (iii) also hold "mutatis-mutandis" for stabilizability with the modification that the perturbed system is stabilizable in Property (i) if and only if it has no input-decoupling zeros in $C_{+0}$.

**(v)** Properties (i)-(iii) also hold for observability by replacing $C^{q_A+q_B+2} \to C^{q_A+q_C+2}$,

$$Q\left(z^{(A)}, z^{(B)}\right) \to Q'\left(z^{(A)}, z^{(B)}\right) := \left(-\sum_{i=0}^{q_A}\sum_{j=1}^{n_A} z_i^{(A)*} D_{ij}^{(A)*} \otimes E^{(A)} \ \vdots \ \sum_{i=0}^{q_C}\sum_{j=1}^{n_C} z_i^{(C)} D_{ij}^{(C)*} \otimes E^{(B)}\right)$$

$\left(-\tilde{A}(z^{(A)}): \tilde{B}(z^{(B)})\right) \to \left(-\tilde{A}(z^{(A)})^* : \tilde{C}(z^{(C)})^*\right)$ and the modification that the perturbed system is observable in Property (i) if and only if it has no output-decoupling zeros in $C$. Properties (i) and (iii) also hold "mutatis-mutandis" for detectability with the changes $C^{q_A+q_B+2} \to C^{q_A+q_C+2}$, $Q\left(z^{(A)}, z^{(B)}\right) \to Q'\left(z^{(A)}, z^{(B)}\right)$, $\left(-\tilde{A}(z^{(A)}): \tilde{B}(z^{(B)})\right) \to \left(-\tilde{A}(z^{(A)})^* : \tilde{C}(z^{(C)})^*\right)$ and $C \to C_{+0}$ in Property (i).

**Proof**: **(i)** Since the nominal system (1)-(2) is controllable in $C_{C0}$,

$$0 < \varepsilon_{C0} \le \sup\left(\left\|\hat{Z}_{C0}\left(s, z^{(A)}, z^{(B)}\right)\right\|_2 : s \in C, \left(z^{(A)}, z^{(B)}\right) \in C_{C0}\right)$$

$$= sup\left(\lambda_{max}\left(\hat{Z}_{C0}\left(s, z^{(A)}, z^{(B)}\right)\right): s \in \overline{C}_{\infty}^{0}, \left(z^{(A)}, z^{(B)}\right) \in C_{C0}^{Fr}\right) \leq \delta_{C0} < \infty \quad (20)$$

where $\hat{Z}_{C0}\left(s, z^{(A)}, z^{(B)}\right) := Z_{C0}\left(s, z^{(A)}, z^{(B)}\right) Z_{C0}^{*}\left(s, z^{(A)}, z^{(B)}\right)$ and $C_{\infty}$ is a circumference of closure (i.e. a circle) of infinity radius in the complex plane centered at the origin. Inequalities (20) hold since the maximum eigenvalue is bounded positive real and equalizes the $\ell_2$-norm because $\hat{Z}_{C0}\left(s, z^{(A)}, z^{(B)}\right)$ is (at least) positive semidefinite in $C \times C_{C0}$ and since the maximum eigenvalue, as being a continuous function on its definition domain $C \times C_{C0}$, reaches its maximum on the boundary of such a domain. Equation (20) implies that

$$0 < \delta_{C0}^{-1} \leq sup\left(\left\|\hat{Z}_{C0}^{-1}\left(s, z^{(A)}, z^{(B)}\right)\right\|_2 : s \in C_{-\infty} \cup C_{+\infty}, \left(z^{(A)}, z^{(B)}\right) \in C_{C0}^{Fr}\right) \leq \varepsilon_{C0}^{-1} < \infty$$

where $C_{-\infty} = (C \setminus C_{+\infty}) \cup \{s = i\omega : \omega \in R\}$, with $C_{+\infty} := \{s \in (C_{\infty} \cap C_{0+})\}$, is a semi circumference of infinity radius centered at the origin of the complex plane. Since $\hat{Z}_{C0}^{-1}\left(s, z^{(A)}, z^{(B)}\right)$ is a strictly proper rational matrix function in s, $\forall\left(z^{(A)}, z^{(B)}\right) \in C_{C0}$, the above supremum on $s \in C_{-\infty} \cup C_{+\infty}$ equalizes the supremum on the imaginary axis of the complex plane. Then, the above constraint becomes identical to

$$0 < \delta_{C0}^{-1} \leq sup\left(\left\|\hat{Z}_{C0}^{-1}\left(i\omega, z^{(A)}, z^{(B)}\right)\right\|_2 : \omega \in R, \left(z^{(A)}, z^{(B)}\right) \in C_{C0}^{Fr}\right) \leq \varepsilon_{C0}^{-1} < \infty \quad (21)$$

This follows for the maximum modulus principle by calculating the maximum of an analytical function on the boundary of $C_{\infty} \times C_{C0}$ which is $\{i\omega : \omega \in R\} \times C_{C0}^{Fr}$. Since the nominal system (1)-(2) is controllable in $C_{C0}$, so that $\hat{Z}_{C0}\left(s, z^{(A)}, z^{(B)}\right)$ is positive definite in $C \times C_{C0}$, then $\hat{Z}_{C}\left(s, z^{(A)}, z^{(B)}\right)$ is positive definite in $C \times C_C \subset C \times C_{C0}$ from (19)-(21) with the domain $C_C$ defined such that $1 > \varepsilon_{C0}^{-1}(2\delta_{C0} + \delta)\delta$ and $sup\left(\left\|\left(-z^{(A)}I_n \vdots z^{(B)}I_n\right) \Delta^{(AB)}\right\|_2 : \left(z^{(A)}, z^{(B)}\right) \in C_C\right) < \delta$ for each $C_C \subset C_{C0}$ from Banach´s Perturbation Lemma, [40]. Thus, $C_C$ exists defined by $\cup Dom\left(\Delta^{(AB)}\right)$ of all the parametrical multi perturbations $\Delta^{(AB)} \in P_{\Delta AB} \cap C^{2n \times (n+m)}$ with:

$$P_{\Delta AB} := \left\{\Delta^{(AB)} : Dom\left(\Delta^{(AB)}\right) \to C^{2n \times (n+m)} : sup\left(-z^{(A)}I_n \vdots z^{(B)}I_n\right)\right\| \left\|\Delta^{(AB)}\right\|_2$$
$$\left(:\left(z^{(A)}, z^{(B)}\right) \in C_C\right) < \delta := \sqrt{\delta_{C0}^2 + \varepsilon_{C0}} - \delta_{C0}\right\}$$

where $P_{\Delta AB}$ depends on $C_C$ and is formed by all the parametrical perturbations $\Delta^{(AB)} \in P_{\Delta AB}$ which satisfy:

$$\bar{\sigma}\left(\Delta^{(AB)}\right) := \left\|\Delta^{(AB)}\right\|_2 < \delta / \sqrt{\sup\left(\left\|z^{(A)}\right\|_2^2 + \left\|z^{(B)}\right\|_2^2 : \left(z^{(A)} \vdots z^{(B)}\right) \in C_c\right)}$$

Thus, the system is controllable within prefixed bounded open domain $C_C \subset C_{C0}$ for all parametrical multi perturbations in the set:

$$P_A \times P_B := \left\{\left(\tilde{A}(z^{(A)}) \vdots \tilde{B}(z^{(B)})\right) \equiv \left(z^{(A)} I_n \vdots z^{(B)} I_n\right) \Delta^{(AB)} : \left(z^{(A)}, z^{(B)}\right) \in C_C, \Delta^{(AB)} \in P_{\Delta AB}\right\}$$

**(ii)** Parametrical perturbations $\left(\tilde{A}(z^{(A)}) \vdots \tilde{B}(z^{(B)})\right)$ with the structure $\left(\sum_{j=1}^{n_A} D_{0j}^{(A)} \Delta_{0j}^{(A)} E^{(A)} \vdots \sum_{j=1}^{n_A} D_{0j}^{(B)} \Delta_{0j}^{(B)} E^{(B)}\right)$ lie in $P_A \times P_B$ from (7). Note that the linear algebraic equation $Q\left(z^{(A)}, z^{(B)}\right) x_Q = v_Q$ has at least a solution for $v_Q$ being of the given form and

$$x_Q := \left(vec^T\left(D_{ij}^{(A)} : (i, j) \in \bar{q}_A^0 \times \bar{n}_A\right), vec^T\left(D_{ij}^{(A)} : (i, j) \in \bar{q}_B^0 \times \bar{n}_B\right)\right)^T$$

Direct calculations in (19) yield that $\det \hat{Z}_C\left(s, z^{(A)}, z^{(B)}\right) = 0$, since $\operatorname{rank} Z_C\left(s, z^{(A)}, z^{(B)}\right) = \operatorname{rank}\left(Z_{C0}\left(s, z^{(A)}, z^{(B)}\right) + \left(-\tilde{A}(z^{(A)}) \vdots \tilde{B}(z^{(B)})\right)\right) < n$ for some $\left(s, \left(z^{(A)}, z^{(A)}\right)\right) \in C \times C_C$ since there exists a nonzero $v_Q$ such that

$$\left(Z_{C0}\left(s, z^{(A)}, z^{(B)}\right) \otimes I_n + Q\left(z^{(A)}, z^{(B)}\right) \otimes x_Q^*\right) v_Q = 0$$

and the proof of Property (ii) is complete.

**(iii)** Since the nominal system is controllable, $0 < \varepsilon_1 \leq \det \hat{Z}_{C0}\left(s, z^{(A)}, z^{(B)}\right) \leq \varepsilon_2 < \infty$ so that $0 < \varepsilon_2^{-1} \leq \det \hat{Z}_{C0}\left(s, z^{(A)}, z^{(B)}\right) \leq \varepsilon_1^{-1} < \infty$, $\forall\left(s, z^{(A)}, z^{(B)}\right) \in C \times C^{q_A + q_B + 2}$. Then, from (19), it follows that:

**1.** There are infinitely many parametrical multi perturbations $\left(-\tilde{A}(z^{(A)}) \vdots \tilde{B}(z^{(B)})\right)$, defined by a map $d : C^{q_A + q_B + 2} \to P_A \times P_B \subset C^{n \times (n+m)}$ satisfying (17)-(18) and fulfilling that $\operatorname{rank}\left(-z^{(A)} I_n \vdots z^{(B)} I_n\right) \Delta^{(AB)} = n$. Note that such parametrical perturbations always exist from (17)-(18) since $E^{(A)}$ and $E^{(B)}$ are full rank.

2. There exist continuous functions $g_1 : C^{q_A+q_B+2} \times P_A \times P_B \to C$ and $g_2 : C^{q_A+q_B+2} \times P_A \times P_B \to C_{+0}$ from (19) satisfying:

$$g_i\left(z^{(A)}, z^{(B)}, \Delta^{(AB)}\right) = O\left(\left\|\Delta^{(AB)}\right\|^i\right) \text{ for } i = 1, 2; \quad \forall \left(z^{(A)}, z^{(B)}\right) \in C^{q_A+q_B+2}$$

$$g_2\left(z^{(A)}, z^{(B)}, \Delta^{(AB)}\right) > 0 \text{ if } \left(z^{(A)}, z^{(B)}\right) \neq 0.$$

Then,

$$\det \hat{Z}_C\left(s, z^{(A)}, z^{(B)}\right) \geq \varepsilon_2^{-1}\left(1 + \text{trace}\left(\varepsilon_3 \varepsilon_1^{-1}\right) + \left|o\left(\varepsilon_3 \varepsilon_1^{-1}\right)\right|\right)$$

$$\times \left\|\left(z^{(A)}, z^{(B)}\right)\right\|\left(\left\|\left(z^{(A)}, z^{(B)}\right)\right\| g_2\left(z^{(A)}, z^{(B)}\right) - \left|g_1\left(z^{(A)}, z^{(B)}\right)\right|\right) \to \infty$$

for some $\varepsilon_3 \in R$ as $C^{q_A+q_B+2} \ni \left(z^{(A)}, z^{(B)}\right) \to \infty$ (the infinity point in $C^{q_A+q_B+2}$) and then the perturbed system (1)-(2) is not controllable (see Remark 3.16) what is a contradiction. The proof of Property (iii) is complete.

**(iv)-(v)** Their proofs are similar to those of properties (i)-(iii) with the given modifications and the replacements $\hat{Z}_{CO}\left(s, z^{(A)}, z^{(C)}\right) \to \hat{Z}_{OO}\left(s, z^{(A)}, z^{(C)}\right)$, $\hat{Z}_C\left(s, z^{(A)}, z^{(C)}\right) \to \hat{Z}_O\left(s, z^{(A)}, z^{(C)}\right)$. □

**Remarks 4.2**. The extension of Theorem 4.1 to output controllability is immediate by invoking Definition 3.3 and Theorem 3.11 (iv). Also, note that if the state-space realization of the nominal system (1)-(2) is minimal in a certain domain, then Theorem 4.1 provides testable conditions to guarantee that the realization is maintained minimal for a set of parametrical perturbations (17)-(18) in a certain domain included in the above one if controllability and observability hold jointly in such a domain. □

**Remarks 4.3**. The various properties might also be investigated by the nominal system defined at some $z_0 = \left(z_0^{(A)}, z_0^{(B)}, z_0^{(C)}, z_0^{(D)}\right) \in C^q$ where the corresponding property holds. For this purpose, the replacements below are used to apply Theorem 4.1:

$$\tilde{A}\left(z^{(A)}\right) \to \overline{\tilde{A}}\left(z^{(A)}\right) := \tilde{A}\left(z^{(A)}\right) + A\left(z^{(A)}\right) - A\left(z_0^{(A)}\right), \quad \tilde{B}\left(z^{(B)}\right) \to \overline{\tilde{B}}\left(z^{(B)}\right) := \tilde{B}\left(z^{(B)}\right) + B\left(z^{(B)}\right) - B\left(z_0^{(B)}\right)$$

$$\tilde{C}\left(z^{(C)}\right) \to \overline{\tilde{C}}\left(z^{(C)}\right) := \tilde{C}\left(z^{(C)}\right) + C\left(z^{(C)}\right) - C\left(z_0^{(C)}\right), \quad \tilde{D}\left(z^{(D)}\right) \to \overline{\tilde{D}}\left(z^{(D)}\right) := \tilde{D}\left(z^{(D)}\right) + D\left(z^{(D)}\right) - D\left(z_0^{(D)}\right)$$

by redefining the parametrical perturbations in certain domains of $C^q$ containing $z$ and $z_0$ and $z = \left(z^{(A)}, z^{(B)}, z^{(C)}, z^{(D)}\right)$ as $\overline{\tilde{A}}\left(z^{(A)}\right)$, $\overline{\tilde{B}}\left(z^{(B)}\right)$, $\overline{\tilde{C}}\left(z^{(C)}\right)$ and $\overline{\tilde{D}}\left(z^{(D)}\right)$. Then, domains where the studied property is kept from the nominal system may be obtained in this way. □

**Remark 4.4**. If the system (1)-(2) is parameterized by real or complex Hermitian matrices then (21) may be simplified to the use of $sup\left(\left\|\hat{Z}_{C0}^{-1}\left(i\omega, z^{(A)}, z^{(B)}\right)\right\|_2 : \omega \in \mathbf{R}_{+0}, \left(z^{(A)}, z^{(B)}\right) \in C_{C0}^{Fr}\right)$ from symmetry considerations.

Observability/Detectability can be investigated by replacing $\hat{Z}_{C0}\left(i\omega, z^{(A)}, z^{(B)}\right)$ with $\hat{Z}_{O0}\left(i\omega, z^{(A)}, z^{(C)}\right)$. In the common case that $D=0$ in (2), i.e. there is no input-output interconnection gain the formula (21) may be reformulated in a simplified way by taking into account that

$$sup\left(\left\|\hat{Z}_{C0}^{-1}\left(s, z^{(A)}, z^{(B)}\right)\right\|_2 : s \in C_{-\infty}, \left(z^{(A)}, z^{(B)}\right) \in C_{C0}^{Fr}\right)$$

$$= sup\left(\left\|\hat{Z}_{C0}^{-1}\left(i\omega, z^{(A)}, z^{(B)}\right)\right\|_2 : \omega \in \mathbf{R}, \left(z^{(A)}, z^{(B)}\right) \in C_{C0}^{Fr}\right)$$

were $C_{+\infty} := \{s \in (C_\infty \cap C_+) \cup \{s = i\omega : \omega \in \mathbf{R}\}\}$ and $C_{-\infty} = (C \setminus C_{+\infty}) \cup \{s = i\omega : \omega \in \mathbf{R}\}$ are semi circumferences of infinity radius since $\hat{Z}_{C0}^{-1}\left(s, z^{(A)}, z^{(B)}\right)$ is a strictly proper rational matrix function in $s$, $\forall \left(z^{(A)}, z^{(B)}\right) \in C_{C0}$. The above identity follows for the maximum modulus principle. If, furthermore, the system (1)-(2) is parameterized by real or complex Hermitian matrices then the above formulas equalize to $= sup\left(\left\|\hat{Z}_{C0}^{-1}\left(i\omega, z^{(A)}, z^{(B)}\right)\right\|_2 : \omega \in \mathbf{R}_{+0}, \left(z^{(A)}, z^{(B)}\right) \in C_{C0}^{Fr}\right)$ by symmetry considerations.

**5. Application to time-delay systems with point internal and external delays**

Now, consider the following extension of the dynamic system (1)-(2) including delays:

$$\dot{x}(z, t) = \left(A\left(z^{(A)}\right) + \tilde{A}\left(z^{(A)}\right)\right)x(z, t) + \left(B\left(z^{(B)}\right) + \tilde{B}\left(z^{(B)}\right)\right)u(t)$$
$$+ \sum_{j=1}^{\eta}\left(A_j\left(z^{(Ad)}\right) + \tilde{A}_j\left(z^{(Ad)}\right)\right)x(z, t-h_j) + \sum_{j=1}^{\kappa}\left(B_j\left(z^{(Bd)}\right) + \tilde{B}_j\left(z^{(Bd)}\right)\right)u(z, t-h'_j)$$

(22)

$$y(z, t) = \left(C\left(z^{(C)}\right) + \tilde{C}\left(z^{(C)}\right)\right)x(z, t) + \left(D\left(z^{(D)}\right) + \tilde{D}\left(z^{(C)}\right)\right)u(t) \qquad (23)$$

fully described by (1)-(2), subject to parametrical multi perturbations (17)-(18), $\eta$ internal (i.e. in the state) pair-wise distinct point delays $h_j \in [0, \bar{h}_j) \subset \mathbf{R}_{+0}$, $\forall j \in \bar{\eta}$ and $\kappa$ external (i.e. in the input) pair-wise distinct point delays $h'_j \in \left[0, \bar{h}'_j\right) \subset \mathbf{R}_{+0}$, $\forall j \in \bar{\kappa}$. The intervals $\left[0, \bar{h}_{(.)}\right)$ and $\left[0, \bar{h}'_{(.)}\right)$ are the admissibility domains of the corresponding internal and external delays, respectively. The external delays could be considered to act on the output instead of on the input with no loss in generality. The initial conditions are defined by any bounded piecewise absolutely continuous vector function with eventual isolated discontinuities $\phi: \left[-\bar{h}, 0\right] \to \mathbf{C}^n$ where $\bar{h} := \underset{1 \le j \le \eta}{Max}(h_j)$. If $u \in PC^{(0)}\left(\mathbf{C}^q \times \mathbf{R}_{+0}, U\right)$ then a unique solution exists on $\mathbf{R}_+$ for each given $\phi: \left[-\bar{h}, 0\right] \to \mathbf{C}^n$, [8], [10], [12]. To set an appropriate framework related to that of the preceding sections, consider $\hat{q}$-tuples $z$ and $\hat{z}_\varphi$ in $\mathbf{C}^{\hat{q}}$ with

$$\hat{q} := q + q_{Ad} + q_{Bd} = q_A + q_B + q_C + q_D + 4 + q_{Ad} + q_{Bd}:$$

$$z := \left(z^{(A)}, z^{(B)}, z^{(C)}, z^{(D)}, z^{(Ad)}, z^{(Bd)}\right); \quad \hat{z}_{\rho\varphi} := \left(z^{(A)}, z^{(B)}, z^{(C)}, z^{(D)}, \hat{z}^{(Ad)}_{\rho\varphi}, \hat{z}^{(Bd)}_{\rho'\varphi'}\right) \quad (24)$$

where

$$\hat{z}^{(Ad)}_{\rho\varphi} := \left(z^{(Ad)}_1 \rho_{1 \prec \varphi_1}, z^{(Ad)}_1 \rho_{2 \prec \varphi_2}, \cdots, z^{(Ad)}_1 \rho_{\eta \prec \varphi_\eta}, \cdots, z^{(Ad)}_{q_{Ad}} \rho_{1 \prec \varphi_1}, z^{(Ad)}_{q_{Ad}} \rho_{2 \prec \varphi_2}, \cdots, z^{(Ad)}_{q_{Ad}} \rho_{\eta \prec \varphi_\eta}\right)$$

$$(25)$$

$$\hat{z}^{(Bd)}_{\rho'\varphi'} := \left(z^{(Ad)}_1 \rho'_{1 \prec \varphi'_1}, z^{(Bd)}_1 \rho'_{2 \prec \varphi'_2}, \cdots, z^{(Bd)}_1 \rho'_{\kappa \prec \varphi'_\kappa}, \cdots, z^{(Bd)}_{q_{Bd}} \rho'_{1 \prec \varphi'_1}, z^{(Bd)}_{q_{Bd}} \rho'_{2 \prec \varphi'_2}, \cdots, z^{(Bd)}_{q_{Bd}} \rho'_{\kappa \prec \varphi'_\kappa}\right)$$

$$(26)$$

and $\hat{q}$-tuples being complex vector functions from $\mathbf{C} \times \mathbf{C}^{\hat{q}}$ to $\mathbf{C}^{\hat{q}}$

$$\bar{z}_{hh'}(s) := \left(z^{(A)}, z^{(B)}, z^{(C)}, z^{(D)}, \bar{z}^{(Ad)}_h(s), \bar{z}^{(Bd)}_{h'}(s)\right) \quad (27)$$

where $\mathbf{h} := \left(h_1, h_2, ..., h_\eta\right)$ and $\mathbf{h}' := \left(h'_1, h'_2, ..., h'_\kappa\right)$ are tuples formed with the sets of internal and external delays, respectively, and

$$\bar{z}^{(Ad)}_h(s) := \left(z^{(Ad)}_1 e^{-h_1 s}, z^{(Ad)}_1 e^{-h_2 s}, \cdots, z^{(Ad)}_1 e^{-h_\eta s}, \cdots, z^{(Ad)}_{q_{Ad}} e^{-h_1 s}, z^{(Ad)}_{q_{Ad}} e^{-h_2 s}, \cdots, z^{(Ad)}_{q_{Ad}} e^{-h_\eta s}\right)$$

$$(28)$$

$$\bar{z}^{(Bd)}_{h'}(s) := \left(z^{(Ad)}_1 e^{-h'_1 s}, z^{(Bd)}_1 e^{-h'_2 s}, \cdots, z^{(Bd)}_1 e^{-h'_\kappa s}, \cdots, z^{(Bd)}_{q_{Bd}} e^{-h'_1 s}, z^{(Bd)}_{q_{Bd}} e^{-h'_2 s}, \cdots, z^{(Bd)}_{q_{Bd}} e^{-h'_\kappa s}\right)$$

$$(29)$$

$\forall \varphi_j, \varphi'_\ell \in [0, 2\pi), \forall \rho_j, \rho'_\ell \in \mathbf{R}_+$, $\left(j \in \bar{\eta}, \ell \in \bar{\kappa}\right)$ where $\gamma_{\prec \theta} := \gamma(\cos\theta + i\sin\theta)$ is a circumference of radius $\gamma$ centered at zero in the complex plane. The dynamic system (22)-(23) may be equivalently described through an algebraic linear system by taking Laplace transforms with zero initial conditions as follows:

$$\left(sI_n - A\left(z^{(A)}\right) + \tilde{A}\left(z^{(A)}\right) - \sum_{j=1}^{\eta}\left(A_j\left(z^{(Ad)}\right) + \tilde{A}_j\left(z^{(Ad)}\right)\right)e^{-h_j s}\right)\hat{x}(z,s)$$

$$-\left(B\left(z^{(B)}\right) + \tilde{B}\left(z^{(B)}\right) - \sum_{j=1}^{K}\left(A_j\left(z^{(Ad)}\right) + \tilde{A}_j\left(z^{(Ad)}\right)\right)e^{-h'_j s}\right)\hat{u}(s) = 0 \quad (30)$$

$$\hat{y}(z,s) = \left(C\left(z^{(C)}\right) + \tilde{C}\left(z^{(C)}\right)\right)\hat{x}(z,s) + \left(D\left(z^{(D)}\right) + \tilde{D}\left(z^{(C)}\right)\right)\hat{u}(s) \quad (31)$$

Note by direct inspection that there exist surjective mapping of the set of tuples (28) to the set of tuples (26) and from the set of tuples (29) to the set of tuples (27) by considering them as functions from $C$ to $[0,\bar{h}_1) \times ... \times [0,\bar{h}_\eta) \times R^{2\eta} \times R_{+0} \times [0,2\pi)$ and to $[0,\bar{h}'_1) \times ... \times [0,\bar{h}'_\eta) \times R^{2\kappa} \times R_{+0} \times [0,2\pi)$, respectively, by associating

$$s = \sigma + i\omega \to \rho_k = e^{-h_k \sigma}, \varphi_k = h_k \omega; \quad \forall k \in \bar{\eta} \quad \text{and} \quad s = \sigma + i\omega \to \rho'_k = e^{-h'_k \sigma}, \varphi'_k = h'_k \omega; \quad \forall k \in \bar{\kappa},$$

respectively. However, those mappings are not one-to-one, in general, for all the admissible sets of delays, since the inverse maps:

$$[0,\bar{h}_1) \times ... \times [0,\bar{h}_\eta) \times R^{2\eta} \times R_{+0} \times [0,2\pi) \subset R^{\eta+2} \to C \text{ and } [0,\bar{h}'_1) \times ... \times [0,\bar{h}'_\kappa) \times R^{2\kappa} \times R_{+0} \times [0,2\pi) \subset R^{\kappa+2} \to C$$

do not have the same definition domain as the respective ranges of the original mappings. The first inverse mapping only exists if and only if $\frac{\ln \rho_k}{h_k}$ are identical and real, $\forall k \in \bar{\eta}$ and also if $\frac{\varphi_k}{h_k}$ are identical, real and belong to $[0,2\pi), \forall k \in \bar{\eta}$. The second inverse mapping only exists if and only if $\frac{\ln \rho'_k}{h'_k} = K \in R$, $\frac{\varphi'_k}{h'_k} = K' \in [0,2\pi)$, $\forall k \in \bar{\kappa}$. The simultaneous existence of both inverse mappings require the fulfillment of joint constraints:

$$\frac{\ln \rho_j}{h_j} = \frac{\ln \rho'_k}{h'_k} = K \in R, \frac{\varphi_j}{h_j} = \frac{\varphi'_k}{h'_k} = K' \in [0,2\pi); \forall j \in \bar{\eta}, \forall k \in \bar{\kappa} \quad (32)$$

**Remark 5.1**. The properties of controllability, observability, stabilizability and detectability of the system (30)-(31) in a domain may be directly tested by extending directly Theorem 3.11, Corollaries 3.12 and 3.14 and Theorem 4.1 under the subsequent guidelines. It turns out that if the constraints (32) are not invoked, only sufficient conditions may be obtained by extending the results of the preceding sections for each tested property by considering the tuples (24)-(26) in the tests. If, in addition, the constraints (32), associated with (27)-(29), are required for given sets of internal and external delays then necessary and sufficient conditions may be obtained by extending such tests from the delay-free case. If the tests fail only for certain sets of $\hat{q}$-tuples (24)-(26) in some given domain, which do not fulfill (32) for given sets of delays, then the system fulfills the tested property for that set of delays. The property is also lost for the sets of delays which do not fulfill the test for the tuples (24)-(26) which have a solution under the constraints (32). If the test does not fail for any tuple (24)–(26) in some domain then the system fulfills the tested property independent of the delays in such a domain. In summary, replace $q \to \hat{q}$, $z^{(A)} \to \left(z^{(A)}, \hat{z}_{\rho\varphi}^{(Ad)}\right)$, $z^{(B)} \to \left(z^{(B)}, \hat{z}_{\rho'\varphi'}^{(Bd)}\right)$ and extend Definitions 3.1-3.4 to the system (30)-(31) to then generalize the various results in Theorem 3.11, Corollaries 3.12 and 3.14 and Theorem 4.1 to the system (30)–(31) subject to delays. Then,

**1.** If any investigated property (namely, controllability, observability, stabilizability or detectability) holds for all z (defined in (24)) in a domain then the property holds within such a domain independent of the delays, i.e. for all sets of $\eta$ internal delays and $\kappa$ external delays ranging from zero to infinity.

**2.** Assume that two sets of internal and external delays are given and assume also that the investigated property holds for all z except at isolated points in a certain domain. Then, if some of the constraints (32) fails for the given sets of delays for all those all points in the domain then the system possesses the investigated property for the given sets of delays. If no sets of delays are specified, then the investigated property except holds for all delays in the admissibility domain except for those where some of the joint constraints (32) fails.

**3.** If the investigated property fails at some z only for sets of delays which fulfill (32) then the system maintains such a property for all delays in their admissibility domains except for those which fulfill (32).

The following technical result is useful for testing the controllability of the time-delay nominal system. The perturbed system is guaranteed to maintain controllability if the nominal one is controllable by incorporating extended sufficiency-type conditions from "ad –hoc" extended versions of Theorem 3.11, Corollaries 3.12 and 3.14 and Theorem 4.1 (see Remark 5.1). For the remaining properties of observability, stabilizability and detectability, the result is extendable "mutatis-mutandis" with the corresponding changes.

**Theorem 5.2** Assume that $C_{\alpha 0}^{(A,B,d)} := C_{\alpha 0}^{(A)} \times C_{\alpha 0}^{(B)} \times C_{\alpha 0}^{(Ad)} \times C_{\alpha 0}^{(Bd)} \subset C^{\hat{q}_C}$ is bounded where $\hat{q}_C := q_A + q_B + q_{Ad} + q_{Bd} + 2$. Then, the following properties hold:

**(i)** The nominal system with delays (22)-(23) is controllable independent of the delays on $C_{\alpha 0}^{(A,B,d)}$ if

$$det\left(\hat{Z}_{C0}(\hat{\sigma}, z)\right) \geq \underline{\delta} > 0, \forall (\hat{\sigma}, z) \in \mathbf{R}^{\eta+\kappa+1} \times [0, 2\pi)^{\eta+\kappa+1} \times C_{\alpha 0}^{(A,B)}$$

where

$$\hat{\sigma} := \left(\sigma, \rho_1, \rho_2, ..., \rho_\eta, \rho'_1, \rho'_2 ..., \rho'_\kappa, \omega, \varphi_1, \varphi_2, ..., \varphi_\eta, \varphi'_1, \varphi'_2 ..., \varphi'_\kappa\right) \in \mathbf{R}^{\eta+\kappa+1} \times [0, 2\pi)^{\eta+\kappa+1}$$

which satisfies all the constraints $\sigma = -\frac{\ln \rho_i}{h_i} = -\frac{\ln \rho'_j}{h'_j}$, $\omega = \frac{\arccos \varphi_i}{h_i} = \frac{\arccos \varphi'_j}{h'_j}$, $\forall (i,j) \in \overline{\eta} \times \overline{k}$ for any two sets of distinct nonnegative real numbers $\{h_1, h_2, ..., h_\eta\}$, $\{h'_1, h'_2, ..., h'_\kappa\}$.

**(ii)** The nominal system with delays (22)-(23) is controllable on $C_{\alpha 0}^{(A,B,d)}$ for a set of $\eta$ positive distinct internal delays $\{h_1, h_2, ..., h_\eta\}$ and a set of $\kappa$ positive distinct external delays $\{h'_1, h'_2, ..., h'_\kappa\}$ provided that

$$\det(\hat{Z}_{C0}(\hat{\sigma}, z)) \geq \underline{\delta} > 0, \forall (\hat{\sigma}, z) \in \mathbf{R}^{\eta+\kappa+1} \times [0, 2\pi)^{\eta+\kappa+1} \times C_{\alpha 0}^{(A,B,d)}$$

where

$$\hat{\sigma} := (\sigma, \rho_1, \rho_2, ..., \rho_\eta, \rho'_1, \rho'_2, ..., \rho'_\kappa, \omega, \varphi_1, \varphi_2, ..., \varphi_\eta, \varphi'_1, \varphi'_2, ..., \varphi'_\kappa) \in \mathbf{R}^{\eta+\kappa+1} \times [0, 2\pi)^{\eta+\kappa+1}$$

satisfies $\sigma = -\frac{\ln \rho_i}{h_i} = -\frac{\ln \rho'_j}{h'_j}$, $\omega = \frac{\arccos \varphi_i}{h_i} = \frac{\arccos \varphi'_j}{h'_j}$, $\forall (i,j) \in \overline{\eta} \times \overline{k}$.

**(iii)** Finite covers of $C_{\alpha 0}^{(A,B,d)}$ can be constructed such that sufficient-type conditions of controllability of Properties (i)-(ii) of the nominal system (22)-(23) on $C_{\alpha 0}^{(A,B)}$ may be constructed involving a finite numbers of computations.

**Proof: (i)-(ii).** The proofs of (i)-(ii) follow directly from the structure of the delay system (22)-(23) and Remark 5.1

**(iii)** It is organized in a very technical way. Note that the matrix function $\hat{Z}_{C0}(s, z^{(A)}, z^{(B)})$ is Hermitian (and thus normal) by construction even if the various matrix functions defining the system (1)-(2) are not Hermitian. Assume that $C_R$ is a bounded open or closed circle of finite radius R centred at the origin of $\mathbf{C}$ and that $R_R^2$ is a circle of the same radius centred at the origin of $\mathbf{R}^2$ which is open if $C_R$ and closed if $C_R$ is closed. Then, there is a natural mapping from $C_R \times C_{\alpha 0}^{(A,B)}$ to $R_R^2 \times C_{\alpha 0}^{(A,B,d)}$ associating each matrix $\hat{Z}_{C0}(s, z^{(A)}, z^{(B)})$ in $\mathbf{C} \times C_{\alpha 0}^{(A,B,d)}$ to a matrix $\hat{Z}_{C0}(\sigma, \omega, z^{(A)}, z^{(B)})$ in $R_R^2 \times C_{\alpha 0}^{(A,B,d)}$ where $\sigma = \text{Re } s$ and $\omega = \text{Im}(s)$. Both matrices have the same spectrum which consists of a set of $\upsilon \leq (n+m)$ real vector functions $\lambda_i : R_R^2 \times C_{\alpha 0}^{(A,B,d)} \to \mathbf{R}_{0+}$ $(i \in \overline{\upsilon})$ of multiplicities $\upsilon_i$ $(i \in \overline{\upsilon})$ which satisfy $\sum_{i=1}^{\upsilon} \upsilon_i = n+m$. The numbers $\lambda_i (i \in \overline{\upsilon})$, $\upsilon$, $\upsilon_i (i \in \overline{\upsilon})$ depend on each tuple in $C_R \times C_{\alpha 0}^{(A,B)}$, equivalently, on each tuple in $R_R \times C_{\alpha 0}^{(A,B,d)}$. If each element of the spectrum $\lambda_i (i \in \overline{\upsilon})$ of multiplicity $\upsilon_i (i \in \overline{\upsilon})$ is considered as $\upsilon_i$ identical elements, then there is a bijective mapping from a continuous vector function $\lambda_C : C_R \times C_{\alpha 0}^{(A,B,d)} \to R_{+0}^{(n+m)}$ to another one $\lambda_{R^2} : R_R^2 \times C_{\alpha 0}^{(A,B,d)} \to R_{+0}^{(n+m)}$. Since the images are identical, it is not made distinction between $\lambda_C$ and $\lambda_{R^2}$ by using simply the notation $\lambda$ for both vector functions in their definition domains. It is obvious that controllability of the nominal system (1)-(2) over $C_{\alpha 0}^{(A,B,d)}$ holds if and only if $\lambda : R_R^2 \times C_{\alpha 0}^{(A,B,d)} \to R_+^{(n+m)}$. Now, if $C_{\alpha 0}^{(A,B,d)}$ is simply connected, closed and bounded, then can be covered by a finite cover from Heine-Borel covering theorem. First, assume that $C_{\alpha 0}^{(A,B,d)}$ is connected, closed and bounded and $R_R^2$ is closed then a finite cover $Co_{\alpha 0}^{(A,B,d)}$ exists for $R_R^2 \times C_{\alpha 0}^{(A,B,d)}$, since $R_R^2 \times C_{\alpha 0}^{(A,B,d)}$ is bounded and closed. $Co_{\alpha 0}^{(A,B,d)}$ is the union of finite covers for each of its L connected components $R_R^2 \times C_{\alpha 0 \ell}^{(A,B,d)}$ $(\ell \in \overline{L})$ with $L < \infty$ since $C_{\alpha 0}^{(A,B,d)}$ is bounded, Then, $Co_{\alpha 0}^{(A,B,d)} = \bigcup_{\ell=1}^{L} \left( R_R^2 \times C_{\alpha 0 \ell}^{(A,B,d)} \right)$. Subsequently, assume that $C_{\alpha 0}^{(A,B,d)}$ is not closed (for instance, open or semi-open) with identical remaining hypotheses as above. Then, there is an open bounded set $Coo_{\alpha 0}^{(A,B,d)} \supset R_R^2 \times C_{\alpha 0}^{(A,B,d)} = \bigcup_{\ell=1}^{L} \left( R_R^2 \times C_{\alpha 0 \ell}^{(A,B,d)} \right)$. Thus, a denumerable cover exists for $Co_{\alpha 0}^{(A,B,d)}$ from Lindelöff covering theorem but it still exists a finite cover $Co_{\alpha 0}^{(A,B,d)}$ of $R_R^2 \times C_{\alpha 0}^{(A,B)}$ satisfying the set inclusion chain:

$$Co_{\alpha 0}^{(A,B,d)} = \bigcup_{\ell=1}^{L} \left( R_R^2 \times C_{\alpha 0 \ell 0}^{(A,B,d)} \right) \supset c\ell\, Co_{\alpha 0}^{(A,B,d)} = \bigcup_{\ell=1}^{L} c\ell \left( R_R^2 \times C_{\alpha 0 \ell}^{(A,B,d)} \right)$$

$$= \bigcup_{\ell=1}^{L} c\ell \left( R_R^2 \times C_{\alpha 0 \ell}^{(A,B,d)} \times \right) \supset Co_{\alpha 0}^{(A,B,d)} \supset \bigcup_{\ell=1}^{L} \left( R_R^2 \times C_{\alpha 0 \ell}^{(A,B,d)} \right) = R_R^2 \times C_{\alpha 0}^{(A,B,d)}$$

from Heine-Borel covering theorem for $c\ell\, Coo_{\alpha 0}^{(A,B,d)}$ (which is bounded and closed) where $C_{\alpha 0 \ell 0}^{(A,B,d)} \supset C_{\alpha 0 \ell}^{(A,B,d)}$ $(\ell \in \overline{L})$ is a finite collection of open bounded sets. As a result, a finite cover $Co_{\alpha 0}^{(A,B,d)}$ of $R_R^2 \times C_{\alpha 0}^{(A,B)}$ always exists if $C_{\alpha 0}^{(A,B,d)}$ is bounded and connected. The finite cover:

$$\hat{H}\left( C_{\alpha 0}^{(A,B)} \right) := \bigcup_{i_j \in \overline{L}_j, j \in \overline{q_A + q_B}} \hat{H}_{i_1,...,i_{q+2}}$$

for some prefixed small real $\varepsilon_0 \in R_+$ is a finite union of disjoint hyper-rectangles defined by

$$\hat{H}_{i_1,\ldots,i_{q+2}} := \left\{ (\sigma, \omega, z) \in R_R^2 \times C_{\alpha 0}^{(A,B,d)} : (\sigma, \omega) \in \bigcup_{(j,k) \in \overline{L}_\sigma \times \overline{L}_\omega} H_{jk}(z) \right.$$

$$\left. z = (z_1, z_2, \ldots, z_{q_A + q_B}), -R_{AB} + \sum_{\ell=1}^{i_j - 1} \varepsilon_\ell < z_j < -R_{AB} + \sum_{\ell=1}^{i_j} \varepsilon_\ell + \varepsilon_0, \forall j \in \overline{q_A + q_B} \right\}$$

$\forall i_j \in \overline{L}_j$, $\forall j \in \overline{q_A + q_B}$ where $R_{AB}$ is bounded positive large real number, where

$$H_{jk}(z) := \left\{ (\sigma, \omega, z) \in R_R^2 \times C_{\alpha 0}^{(A,B,d)} \right.$$

$$\left. : -R_R + \sum_{\ell=1}^{j-1} \varepsilon_{\sigma\ell} < \sigma < -R_R + \sum_{\ell=1}^{j} \varepsilon_{\sigma\ell} + \varepsilon_0, -R_R + \sum_{\ell=1}^{k-1} \varepsilon_{\omega\ell} < \omega < -R_R + \sum_{\ell=1}^{k} \varepsilon_{\omega\ell} + \varepsilon_0 \right\}$$

$\forall z \in C_{\alpha 0}^{(A,B,d)}$, $\forall j \in \overline{L}_\sigma$, $\forall k \in \overline{L}_\omega$. Then, $\bigcup_{(j,k) \in \overline{L}_\sigma \times \overline{L}_\omega} H_{jk}(z)$ is a finite cover composed of open rectangles over a large open square in $\mathbf{R}^2$ for each $z$ in $C_{\alpha 0}^{(A,B,d)}$. It is quite obvious that, since any characteristic zeros of the nominal system (1)-(2) are at a finite distance from the origin if the characteristic equation has a principal term and the characteristic quasi-polynomial is monic in the Laplace argument s then the stability may be tested for any z over the union of $H_{jk}(z)$, for $j \in \overline{L}_\sigma < \infty$ and over $\hat{H}\left(C_{\alpha 0}^{(A,B,d)}\right)$ for any parameterization in $C_{\alpha 0}^{(A,B,d)}$. If the nominal system is controllable then the open rectangles $H_{jk}(z)$ for any $z = (z^{(A)}, z^{(B)}) \in C_{\alpha 0}^{(A,B,d)}$ are constructed as follows. Define:

$$f_0(\sigma, \omega, z^{(A)}, z^{(B)}) := \det\left(\hat{Z}_{CO}(\sigma, \omega, z^{(A)}, z^{(B)})\right)$$

Since the above function is analytic everywhere in its definition domain then for given strictly increasing real sequences $\{\sigma_i\}_0^\infty, \{\omega_i\}_0^\infty, \{z_{ji}\}_0^\infty$ $(j \in \overline{q_A + q_B})$ and for any real numbers $\sigma \in (\sigma_i, \sigma_{i+1})$, $\omega \in (\omega_i, \omega_{i+1})$, $z_j \in (z_{ji}, z_{j,i+1})$, $\forall j \in \overline{q_A + q_B}$ being the components of $z = (z^{(A)}, z^{(B)})$:

$$\left| f_0(\sigma, \omega, z^{(A)}, z^{(B)}) \right|_i \geq \left| f_0(\sigma_i, \omega_i, z^{(Ai)}, z^{(Bi)}) \right| - (\sigma - \sigma_i) \left| f'_{0\sigma}(\sigma_i, \omega_i, z^{(A)}, z^{(B)}) \right|$$

$$- (\omega - \omega_i) \left| f'_{0\omega}(\sigma_i, \omega_i, z^{(A)}, z^{(B)}) \right| - \sum_{j=1}^{q_A + q_B} (z_j - z_{ji}) \left| f'_{0z_j}(\sigma, \omega, z^{(A)}, z^{(B)}) \right| \geq \delta_i > 0$$

provided that the above strictly increasing sequences are chosen subject to $(\sigma_i, \omega_i, z_{1i}, \ldots, z_{q_A + q_B, i}) \in R_R^2 \times C_{\alpha 0}^{(A,B,d)}$

$\sigma_i < \sigma_{i+1} < \sigma_i + \varepsilon_i$, $\omega_i < \omega_{i+1} < \omega_i + \varepsilon_i$, $z_{ji} < z_{j,i+1} < z_{ji} + \varepsilon_i$, $\forall j \in \overline{q_A + q_B}$

for any real $\varepsilon_i$ fulfilling $0 < \varepsilon_i \leq \dfrac{\delta_i}{(q_A + q_B + 2) M_i}$, where

$$M_i := \max\left( \left| f'_{0\sigma}\left(\sigma, \omega, z^{(A)}, z^{(B)}\right) \right|, \left| f'_{0\omega}\left(\sigma_i, \omega_i, z^{(A)}, z^{(B)}\right) \right|, \left| f'_{0z_j}\left(\sigma, \omega, z^{(A)}, z^{(B)}\right) \right|, j \in \overline{q_A + q_B} \right)$$

which is finite since $R_R^2 \times C_{\alpha 0}^{(A,B,d)}$ is bounded. Then

$$\hat{H}\left(C_{\alpha 0}^{(A,B,d)}\right) := \bigcup_{i_j \in \overline{L}_j, \, j \in \overline{q_A + q_B}} \hat{H}_{i_1, \ldots, i_{q+2}} \supseteq R_R^2 \times C_{\alpha 0}^{(A,B,d)}$$ is a finite cover of

$R_R^2 \times C_{\alpha 0}^{(A,B,d)}$ provided that $\hat{H}_{i_1, \ldots, i_{q+2}}$ are hyper-rectangles defined by

$\sigma_i < \sigma_{i+1} < \sigma_i + \varepsilon_i$, $\omega_i < \omega_{i+1} < \omega_i + \varepsilon_i$, $z_{ji} < z_{j,i+1} < z_{ji} + \varepsilon_i$, $\forall j \in \overline{q_A + q_B}$ and that the real sequence $\{\delta_i\}_0^\infty$ is bounded and positive satisfying $\delta_i \geq \underline{\delta}$ for some real $\underline{\delta} > 0$. As a result, the existence of such a finite cover implies the controllability of the nominal system (1)-(2). If the construction of the cover satisfying the condition $0 < \varepsilon_i \leq \dfrac{\delta_i}{(q_A + q_B + 2) M_i}$ with $\delta_i \geq \underline{\delta} > 0$ is impossible then the nominal system is not controllable.

**Remark 5.3**. Theorem 5.2 combined with Remark 5.1 yields direct sufficiency type conditions for controllability of the system (22)-(23) on some bounded $C_{\alpha}^{(A,B,d)} \subset C_{\alpha 0}^{(A,B,d)}$ independent of or dependent on the delays when subject to parametrical multi-perturbations. Note that as alternative to tests on determinants, tests on the singular values of $Z_{C0}(\hat{\sigma}, z)$, or on the eigenvalues or matrix ranks of $\hat{Z}_{C0}(\hat{\sigma}, z)$, may be applied. Extensions involving input and output decoupling zeros are also direct but specific derivation is omitted by space reasons. On the other hand, the remaining properties may be investigated as well with simple modifications of Theorem 5.2. In particular,

(1) Stabilizability tests follow from Theorem 5.2 for $\sigma \in \mathbf{R}_{+0}$ only. Output controllability tests follow by replacing $Z_{C0}(\hat{\sigma}, z)$ with $Z_{OC0}(\hat{\sigma}, z)$ in Theorem 5.2.

(2) Observability tests on some bounded set $C_{\alpha 0}^{(A,C,d)}$ follow by replacing $C_{\alpha 0}^{(A,B,d)}$ with $C_{\alpha 0}^{(A,C,d)} := C_{\alpha 0}^{(A)} \times C_{\alpha 0}^{(C)} \times C_{\alpha 0}^{(Ad)} \times C_{\alpha 0}^{(Cd)}$ and then using $\hat{Z}_{O0}(\hat{\sigma}, z)$ instead of $\hat{Z}_{C0}(\hat{\sigma}, z)$ in Theorem 5.2.

(3) Detectability tests on some bounded set $C_{\alpha 0}^{(A,C,d)}$ follow by replacing $C_{\alpha 0}^{(A,B,d)}$ with $C_{\alpha 0}^{(A,C,d)}$ and then using $\hat{Z}_{O0}(\hat{\sigma}, z)$ instead of $\hat{Z}_{C0}(\hat{\sigma}, z)$ in Theorem 5.2 for $\sigma \in \mathbf{R}_{+0}$ only.

Theorem 5.2 can also be applied for testing the controllability nominal delay-free system (1)-(2). However, in this simple case the conditions may be investigated by testing a finite number of (in general) complex eigenvalues of $\hat{Z}_{C0}(\hat{\sigma},z)$, or (real) singular values of $Z_{C0}(\hat{\sigma},z)$. In this case, the construction of a finite subcover is easier than that involved in Theorem 5.2 since the above number of eigenvalues/singular values is finite for each point in a bounded set $C_{\alpha 0}^{(A,B)}$ instead of a finite number o functions with infinitely many associated point eigenvalues. □

ACKNOWLEDGMENTS

The author is very grateful to MCYT by its partial support of this work through Grant DPI2006-00714. He is grateful to the Basque Government by their respective partial support of this work via Project DPI 2006-00714, and Research Grants: Research Groups No. IT-269-07 and SAIOTEK 2006/ S-PED06UN10 and S-PE-07UN04. He is also grateful to the referees by their useful comments.